\documentclass[twoside]{article}
\usepackage{amsmath}
\usepackage{amsfonts}
\usepackage{amssymb}
\usepackage{color}

\newtheorem{rem}{Remark}[section]
\newtheorem{defi}{Definition}[section]
\newtheorem{lem}{Lemma}[section]

\newtheorem{pro}{Proposition}[section]
\newtheorem{theo}{Theorem}[section]

\newtheorem{ass}{Assumption}

\usepackage[sc]{mathpazo} 
\usepackage[T1]{fontenc} 
\usepackage{microtype} 

\usepackage[margin=3.5cm,hmarginratio=1:1,top=32mm,columnsep=20pt]{geometry} 
\usepackage{multicol} 
\usepackage[hang, small,labelfont=bf,up,textfont=it,up]{caption} 
\usepackage{booktabs} 
\usepackage{hyperref} 

\usepackage{lettrine} 
\usepackage{paralist} 

\usepackage{abstract} 

\usepackage{titlesec} 
\renewcommand\thesection{\Roman{section}} 
\renewcommand\thesubsection{\Roman{subsection}} 
\titleformat{\section}[block]{\large\scshape\centering}{\thesection.}{1em}{} 
\titleformat{\subsection}[block]{\large}{\thesubsection.}{1em}{} 

\usepackage{fancyhdr} 
\title{\vspace{-15mm}\fontsize{24pt}{10pt}\selectfont\textbf{Long signal change-point detection}} 

\date{}

\usepackage[pdftex]{graphicx}
\usepackage{float}
\usepackage[T1]{fontenc}
\usepackage[english]{babel}
\usepackage[utf8x]{inputenc}
\usepackage{ucs}
\usepackage{amsmath}
\usepackage{amsfonts}
\usepackage{amssymb}
\usepackage{array}
\usepackage{listings}
\usepackage{color}
\usepackage{authblk}

\author[1,2]{G\'erard Biau}
\author[3,4]{Kevin Bleakley}
\author[5]{David M. Mason}
\affil[1]{Sorbonne Universit\'es, UPMC, France}
\affil[2]{Institut Universitaire de France}
\affil[3]{INRIA Saclay, France}
\affil[4]{D\'epartement de Math\'ematiques d'Orsay, France}
\affil[5]{University of Delaware, USA}

\begin{document}

\maketitle

\begin{abstract}
\noindent The detection of change-points in a spatially or time-ordered data
sequence is an important problem in many fields such as genetics and
finance. We derive the asymptotic distribution of a statistic recently
suggested for detecting change-points. Simulation of its estimated limit
distribution leads to a new and computationally efficient change-point
detection algorithm, which can be used on very long signals. We assess the
algorithm experimentally under various conditions.
\end{abstract}


\section{Introduction}

When met with a  data set ordered by time or space, it is often
important to predict when or where something \textquotedblleft
changed\textquotedblright\ as we move temporally or spatially through it. In
biology, for example, changes in an array Comparative Genomic Hybridization
(aCGH) or Chip-Seq data signal as one moves across the genome can represent
an event such as a change in genomic copy number, which is extremely
important in cancer gene detection \cite%
{Picard2005statistical,Shah2007Modeling}. In the financial world, detecting
changes in multivariate time-series data is important for decision-making 
\cite{Talih2005Structural}. Change-point detection can also be used to
detect financial anomalies \cite{bolton2002} and significant changes in a
sequence of images \cite{kim2009}.

Change-point detection analysis is a well-studied field and there are
numerous approaches to the problem. Its extensive literature ranges from
parametric methods using log-likelihood functions \cite%
{carlin1992,Lavielle2006Detection} to nonparametric ones based on
Wilcoxon-type statistics, U-statistics and sequential ranks. The reader is
referred to the monograph \cite{csorgo1997limit} for an in-depth treatment
of these methods.

In change-point modeling it is generally supposed that we are dealing with a
random process evolving in time or space. The aim is to develop a method to
search for a point where possible changes occur in the mean, variance,
distribution, etc. of the process. All in all, this comes down to finding
ways to decide whether a given signal can be considered homogeneous in a
statistical (stochastic) sense.

The present article builds upon an interesting nonparametric change-point
detection method that was recently proposed by Matteson and James \cite%
{matteson2014nonparametric}. It uses U-statistics (see \cite%
{hoeffding1948class}) as the basis of its change-point test. Its interest
lies in its ability to detect quite general types of change in distribution.
Several theoretical results are presented in \cite{matteson2014nonparametric}
to highlight some of the mathematical foundations of their method. These in
turn lead to a simple and useful data-driven statistical test for
change-point detection. The authors then apply this test successfully to
simulated and real-world data.

There are however several weaknesses in \cite{matteson2014nonparametric}
both from theoretical and practical points of view. Certain fundamental
theoretical considerations are incompletely treated, especially the assertion
that a limit distribution exists for the important statistic, upon which the
rest of the approach hangs. On the practical side, the method is
computationally prohibitive for signals of more than a few thousand points,
which is unfortunate because real-world signals can be typically much longer.

Our paper has two main objectives. First, it fills in missing theoretical
results in \cite{matteson2014nonparametric} including a derivation of the
limit distribution of the statistic. This requires the effective application
of large sample theory techniques, which were developed to study degenerate
U-statistics. Second, we provide a method to simulate from an approximate
version of the limit distribution. This leads to a new computationally
efficient strategy for change-point detection that can be run on much longer
signals.

The article is structured as follows. In Section \ref{theoretical} we
provide some context and present the main theoretical results. In Section %
\ref{practical} we show how to approximate the limit distribution of the
statistic, which leads to a new test strategy for change-point detection. We
then show how to extend the method to much longer sequences. Simulations  are provided in Section \ref{applications}. A short
discussion follows in Section \ref{discussion}, and a proof of the paper's
main result is given in Section \ref{maintheorem}. Some important technical
results are detailed in the Appendix.

\section{Theoretical results\label{theoretical}}

\subsection{Measuring differences between multivariate distributions}

Let us first briefly describe the origins of the nonparametric change-point
detection method described in \cite{matteson2014nonparametric}. For random
variables $Y, Z$ taking values in $\mathbb{R}^{d}$, $d\geq 1$, let $\phi
_{Y}$ and $\phi _{Z}$ denote their respective characteristic functions. A
measure of the divergence (or \textquotedblleft difference\textquotedblright
) between the distributions of $Y$ and $Z$ is as follows: 
\begin{equation*}
\mathcal{D}(Y,Z)=\int_{\mathbb{R}}\left\vert \phi _{Y}(t)-\phi
_{Z}(t)\right\vert ^{2}w(t)dt,
\end{equation*}%
where $w(t)$ is an arbitrary positive weight function for which this
integral exists. It turns out that for the specific weight function 
\begin{equation*}
w(t;\beta )=\left( \frac{2\pi ^{1/2}\Gamma (1-\beta /2)}{\beta 2^{\beta
}\Gamma \left( (d+\beta )/2\right) }|t|^{d+\beta }\right) ^{-1},
\end{equation*}%
which depends on a $\beta \in (0,2)$, one can obtain a not immediately obvious but very useful result. 
Let $Y, Y^{\prime }$ be i.i.d. $F_{Y}$ and $Z, Z^{\prime }$ be
i.i.d. $F_{Z}$, with $Y, Y^{\prime }$, $Z$ and $Z^{\prime }$ mutually
independent. Denote by $\left\vert \cdot \right\vert $ the Euclidean norm on 
$\mathbb{R}^{d}$. Then, if 
\begin{equation}
\mathbb{E}(\left\vert Y\right\vert ^{\beta }+\left\vert Z\right\vert ^{\beta
})<\infty ,  \label{betafinite}
\end{equation}%
Theorem 2 of \cite{szekely2005hierarchical} yields that 
\begin{equation}
\mathcal{D}(Y,Z;\beta )=\mathcal{E}(Y,Z;\beta ):=2\mathbb{E}\left\vert
Y-Z\right\vert ^{\beta }-\mathbb{E}\left\vert Y-Y^{\prime }\right\vert
^{\beta }-\mathbb{E}\left\vert Z-Z^{\prime }\right\vert ^{\beta }\geq 0,
\label{E}
\end{equation}%
where we have written $\mathcal{D}(Y,Z;\beta )$ instead of $\mathcal{D}(Y,Z)$
to highlight dependence on $\beta $. Therefore (\ref{betafinite}) implies
that $\mathcal{E}(Y,Z;\beta )\in \lbrack 0,\infty )$. Furthermore, Theorem 2
of \cite{szekely2005hierarchical} says that $\mathcal{E}(Y,Z;\beta )=0$ if
and only if $Y$ and $Z$ have the same distribution. This remarkable result
leads to a simple data-driven divergence measure for distributions. Seen in
the context of hypothesizing a change-point in a signal of independent
observations ${\bf X}=(X_{1},\hdots,X_{n})$ after the $k$-th observation $X_{k}$, we
simply calculate an empirical version of (\ref{E}): 
\begin{align}
\mathcal{E}_{k,n}({\bf X};\beta ) &=
\frac{2}{k(n-k)}\sum_{i=1}^{k}\sum_{j=k+1}^{n}\left\vert
X_{i}-X_{j}\right\vert ^{\beta }-{\binom{k}{2}}^{-1}\sum_{1\leq i<j\leq
k}\left\vert X_{i}-X_{j}\right\vert ^{\beta } \nonumber \\   & \qquad -{\binom{n-k}{2}}%
^{-1}\sum_{1+k\leq i<j\leq n}\left\vert X_{i}-X_{j}\right\vert ^{\beta }.
\label{empdiv}
\end{align}%
Matteson and James \cite{matteson2014nonparametric} state without proof that
under the null hypothesis of $X_{1},\hdots,X_{n}$ being i.i.d. (no
change-points), the sample divergence given in (\ref{empdiv}) scaled by $%
\frac{k(n-k)}{n}$ converges in distribution to a non-degenerate random
variable as long as $\min \{k,n-k\}\rightarrow \infty $. Furthermore, they
state that if there is a change-point between two distinct i.i.d.
distributions after the $k$-th point, the sample divergence scaled by $\frac{%
k(n-k)}{n}$ tends a.s. to infinity as long as $\min \{k,n-k\}\rightarrow
\infty $. These claims clearly point to a useful statistical test for
detecting change-points. However, we cannot find rigorous mathematical
arguments to substantiate them in \cite{matteson2014nonparametric}, nor in
the earlier work \cite{szekely2005hierarchical}.

As this is of fundamental importance to the theoretical and practical
validity of this change-point detection method, we shall show the existence
of the non-degenerate random variable hinted at in \cite%
{matteson2014nonparametric} by deriving its distribution. Our approach
relies on the asymptotic behavior of U-statistic type processes, which were
introduced for the first time for change-point detection in random sequences
in \cite{csorgHo1988invariance}; see also Chapter 2 of the book \cite%
{csorgo1997limit}. We also show that in the presence of a change-point the
correctly-scaled sample divergence indeed tends  to infinity with probability 1.

\subsection{Main result}

Let us first begin in a more general setup. Let $X_{1},\hdots,X_{n}$ be
independent $\mathbb{R}^{d}$-valued random variables. For any symmetric
measurable function $\varphi :\mathbb{R}^{d}\times \mathbb{R}^{d}\rightarrow 
\mathbb{R}$, whenever the indices make sense we define the following four
terms: 
\begin{align*}
V_{k}(\varphi )& :=\sum_{i=1}^{k}\sum_{j=k+1}^{n}\varphi (X_{i},X_{j}), \\
U_{n}\left( \varphi \right) & :=\sum_{1\leq i<j\leq n}\varphi (X_{i},X_{j}),
\\
U_{k}^{(1)}(\varphi )& :=\sum_{1\leq i<j\leq k}\varphi (X_{i},X_{j}), \\
U_{k}^{(2)}(\varphi )& :=\sum_{k+1\leq i<j\leq n}\varphi (X_{i},X_{j}).
\end{align*}%
Otherwise, define the term to be zero; for instance, $U_{1}^{(1)}(\varphi
)=0 $ and $U_{k}^{(2)}(\varphi )=0$ for $k=n-1$ and $n$. Note that in the
context of the change-point algorithm we have in mind, $\varphi(x,y) =
\varphi_{\beta }(x,y):=|x-y|^{\beta }$, $\beta \in (0,2)$, but the following
results are valid for the more general $\varphi $ defined above. Notice also
that the last three terms are U-statistics absent their normalization
constants. Next, let us define 
\begin{equation*}
U_{k,n}(\varphi ):=\frac{2}{k(n-k)}V_{k}(\varphi )-{\binom{k}{2}}%
^{-1}U_{k}^{(1)}(\varphi )-{\binom{n-k}{2}}^{-1}U_{k}^{(2)}(\varphi ).
\end{equation*}
Observe that $U_{k,n}(\varphi )$ is a general version of the empirical
divergence given in (\ref{empdiv}). 
Notice that 
\begin{equation}
V_{k}\left( \varphi \right) =U_{n}\left( \varphi \right) -U_{k}^{\left(
1\right) }\left( \varphi \right) -U_{k}^{\left( 2\right) }\left( \varphi
\right) .  \label{Vk}
\end{equation}%
While $U_{k,n}(\varphi )$ is not a U-statistic, we can use (\ref{Vk}) to
express it as a linear combination of U-statistics. Indeed, we find that 
\begin{equation*}
U_{k,n}(\varphi )=\frac{2(n-1)}{k(n-k)}\left( \frac{U_{n}(\varphi )}{n-1}%
-\left( \frac{U_{k}^{(1)}(\varphi )}{k-1}+\frac{U_{k}^{(2)}(\varphi )}{n-k-1}%
\right) \right) .
\end{equation*}%
Therefore, we now have an expression for $U_{k,n}(\varphi )$ made up of
U-statistics, which will be useful in the following.

Our aim is to use a test based on $U_{k,n}(\varphi )$ for the null
hypothesis $\mathcal{H}_{0}:X_{1},\hdots,X_{n}$ have the same distribution,
versus the alternative hypothesis $\mathcal{H}_{1}$ that there is a
change-point in the sequence $X_{1},\hdots,X_{n}$, i.e., 
\begin{align*}
\mathcal{H}_{1}:& \mbox{ There is a }\gamma \in (0,1)\mbox{ such that }%
\mathbb{P}(X_{1}\leq t)=\cdots =\mathbb{P}(X_{\lfloor n\gamma \rfloor }\leq
t), \\
& \mathbb{P}(X_{\lfloor n\gamma \rfloor +1}\leq t)=\cdots =\mathbb{P}%
(X_{n}\leq t),\,\,\,t\in \mathbb{R}^{d}, \\
& \mbox{and }\mathbb{P}(X_{\lfloor n\gamma \rfloor }\leq t_{0})\neq \mathbb{P%
}(X_{\lfloor n\gamma \rfloor +1}\leq t_{0})\mbox{ for some }t_{0}.
\end{align*}%
For $u$, $v\in \mathbb{R}^{d}$, $u\leq $ $v$ means that each component of $u$
is less than or equal to the corresponding component of $v$. Also note that
for any $z\in \mathbb{R}$, $\lfloor z\rfloor $ stands for its integer part.

Let us now examine the asymptotic properties of $U_{k,n}(\varphi )$. We
shall be using notation, methods and results from Section 5.5.2 of monograph 
\cite{serfling1980} to provide the groundwork. In the following, we shall
denote by $F$ the common (unknown) distribution function of the $X_{i}$
under $\mathcal{H}_{0}$,\thinspace\ $X$ a generic random variable with
distribution function $F$, and $X^{\prime }$ an independent copy of $X$. We
assume that 
\begin{equation}
\mathbb{E}\varphi ^{2}(X,X^{\prime })=\int_{\mathbb{R}^{d}}\int_{\mathbb{R}%
^{d}}\varphi ^{2}(x,y)dF(x)dF(y)<\infty ,  \label{sec}
\end{equation}%
and set $\Theta =\mathbb{E}\varphi (X,X^{\prime })$. We also denote $\varphi
_{1}(x)=\mathbb{E}\varphi (x,X^{\prime })$, and define 
\begin{equation}
h(x,y)=\varphi (x,y)-\varphi _{1}(x)-\varphi _{1}(y),\quad \tilde{h}%
_{2}(x,y)=h(x,y)+\Theta .  \label{H1}
\end{equation}%
With this notation, we see that $\mathbb{E}h(X,X^{\prime })=-\Theta $, and
therefore that $\mathbb{E}\tilde{h}_{2}(X,X^{\prime })=0$. Furthermore, 
\begin{equation}
U_{k,n}(\varphi )=U_{k,n}(h)=U_{k,n}(\tilde{h}_{2}),  \label{eq}
\end{equation}%
since%
\begin{equation*}
\frac{U_{n}(\Theta )}{n-1}-\left( \frac{U_{k}^{(1)}(\Theta )}{k-1}+\frac{%
U_{k}^{(2)}(\Theta )}{n-k-1}\right) =\frac{U_{n}(\psi )}{n-1}-\left( \frac{%
U_{k}^{(1)}(\psi )}{k-1}+\frac{U_{k}^{(2)}(\psi )}{n-k-1}\right) =0,
\end{equation*}%
where $\psi (x,y):=\varphi _{1}(x)+\varphi _{1}(y)$. As in Section 5.5.2 of 
\cite{serfling1980}, we then define the operator $A$ on $L_{2}(\mathbb{R}%
^{d},F)$ by 
\begin{equation}
Ag(x):=\int_{\mathbb{R}^{d}}\tilde{h}_{2}(x,y)g(y)dF(y),\quad x\in \mathbb{R}%
^{d},\,g\in L_{2}(\mathbb{R}^{d},F).  \label{op}
\end{equation}%
Let $\lambda _{i}$, $i\geq 1$, be the eigenvalues of this operator $A$ with
corresponding orthonormal eigenfunctions $\phi _{i}$, $i\geq 1$. Since for
all $x\in \mathbb{R}^{d}$, 
\begin{equation*}
\int_{\mathbb{R}^{d}}\tilde{h}_{2}(x,y)dF(y)=0,
\end{equation*}%
we see with $\phi _{1}:=1$, $A\phi _{1}=0=:\lambda _{1}\phi _{1}.$ Thus $%
( 0,1) =( \lambda _{1},\phi _{1}) $ is an
eigenvalue and normalized eigenfunction pair of the operator $A$. This
implies that for every eigenvalue and normalized eigenfunction pair $\left(
\lambda _{i},\phi _{i}\right) $, $i\geq 2,$ where $\lambda _{i}$ is nonzero, 
\begin{equation*}
\mathbb{E}\left( \phi _{1}(X)\phi _{i}(X)\right) =\mathbb{E}\phi _{i}(X)=0.
\end{equation*}%
Moreover, we have that in $L_{2}(\mathbb{R}^{d}\times \mathbb{R}^{d},F\times F)$, 
\begin{equation*}
\tilde{h}_{2}(x,y)=\lim_{K\rightarrow \infty }\sum_{i=1}^{K}\lambda _{i}\phi
_{i}(x)\phi _{i}(y).
\end{equation*}%
From this we get that%
\begin{equation}
\mathbb{E}\tilde{h}_{2}^{2}(X,X^{\prime })=\sum_{i=1}^{\infty }\lambda
_{i}^{2}.  \label{lam}
\end{equation}%
For further details and theoretical justification of these claims, refer to
Section 5.5.2 of \cite{serfling1980} and both Exercise 44 on pg. 1083 and
Exercise 56 on pg. 1087 of \cite{dunford1988linear}. In fact, we shall
assume further that 
\begin{equation}
\sum_{i=1}^{\infty }|\lambda _{i}|<\infty .   \label{finite}
\end{equation}%
It is crucial for the change-point testing procedure that we shall propose
that the function $\tilde{h}_{2}(x,y)$ defined as in (\ref{finite}) with $%
\varphi (x,y)=\varphi _{\beta }(x,y)=|x-y|^{\beta }$, $\beta \in (0,2)$,
satisfies (\ref{finite}) whenever (\ref{sec}) holds. A proof of this is
given in the Appendix.

Next, for any fixed $\frac{2}{n}\leq t<1-\frac{2}{n}$,\thinspace\ $n\geq 3$,
set 
\begin{align}
\mathbb{Y}_{n}(\tilde{h}_{2},t)& :=\frac{\left( \lfloor nt\rfloor \left(
n-\lfloor nt\rfloor \right) \right) ^{2}}{n^{2}(n-1)}U_{\lfloor nt\rfloor
,n}(\tilde{h}_{2})  \label{correctnor} \\
& =\frac{2\lfloor nt\rfloor \left( n-\lfloor nt\rfloor \right) }{n^{2}}%
\left( \frac{U_{n}(\tilde{h}_{2})}{n-1}-\left( \frac{U_{\lfloor nt\rfloor
}^{(1)}(\tilde{h}_{2})}{\lfloor nt\rfloor -1}+\frac{U_{\lfloor nt\rfloor
}^{(2)}(\tilde{h}_{2})}{n-\lfloor nt\rfloor -1}\right) \right) .  \notag
\end{align}%
We define $U_{0}^{(1)}(\tilde{h}_{2})=0$, $U_{0}^{( 2) }( 
\tilde{h}_{2}) =U_{n}( \tilde{h}_{2}) $, $%
U_{1}^{(1)}(\tilde{h}_{2})/0=0$, and $U_{n-1}^{(2)}(\tilde{h}_{2})/0=0$ ,
which gives%
\begin{equation*}
\mathbb{Y}_{n}(\tilde{h}_{2},t)=0,\text{ for }t\in \left[ 0,\frac{1}{n}%
\right) ,
\end{equation*}%
\begin{equation*}
\mathbb{Y}_{n}(\tilde{h}_{2},t)=\frac{2(n-1)}{n^{2}}\left( \frac{U_{n}(%
\tilde{h}_{2})}{n-1}-\frac{U_{1}^{(2)}(\tilde{h}_{2})}{n-2}\right) ,\text{
for }t\in \left[ \frac{1}{n},\frac{2}{n}\right) ,
\end{equation*}%
\begin{equation*}
\mathbb{Y}_{n}(\tilde{h}_{2},t)=\frac{4(n-2)}{n^{2}}\left( \frac{U_{n}(%
\tilde{h}_{2})}{n-1}-\frac{U_{n-2}^{(1)}(\tilde{h}_{2})}{n-3}-U_{n-2}^{(2)}(%
\tilde{h}_{2})\right) ,\text{ for }t\in \left[ 1-\frac{2}{n},1-\frac{1}{n}%
\right) ,
\end{equation*}%
\begin{equation*}
\mathbb{Y}_{n}(\tilde{h}_{2},t)=\frac{2(n-1)}{n^{2}}\left( \frac{U_{n}(%
\tilde{h}_{2})}{n-1}-\frac{U_{n-1}^{(1)}(\tilde{h}_{2})}{n-2}\right) ,\text{
for }t\in \left[ 1-\frac{1}{n},1\right) ,\text{ and }\mathbb{Y}_{n}(\tilde{h}%
_{2},1)=0.
\end{equation*}%
One can readily check that $\mathbb{Y}_{n}(\tilde{h}_{2},\cdot )\in
D^{1}[0,1]$, the space of bounded measurable real-valued functions defined
on $\left[ 0,1\right] $ that are right-continuous with
left-hand limits. Notice that on account of (\ref{eq}) we can also write $%
\mathbb{Y}_{n}(\tilde{h}_{2},\cdot )=\mathbb{Y}_{n}(\varphi ,\cdot )$, and
we will do so from now on. In the following theorem, $\{\mathbb{B}%
^{(i)}\}_{i\geq 1}$ denotes a sequence of independent standard Brownian
bridges. 
\begin{theo}
\label{theo1} Whenever $X_i$, $i\geq 1$ are i.i.d. $F$ and $\varphi $ satisfies (\ref{sec}) and (\ref{finite}), $%
\mathbb{Y}_{n}(\varphi ,\cdot )$ converges weakly in $D^{1}[0,1]$ to the
tied down mean zero continuous process $\mathbb{Y}$ defined on $[0,1]$ by 
\begin{equation*}
\mathbb{Y}(t):=\sum_{i=1}^{\infty }\lambda _{i}\left( t\left( 1-t\right)
-\left( \mathbb{B}^{\left( i\right) }\left( t\right) \right) ^{2}\right) .
\end{equation*}%
In particular, 
\begin{equation*}
\sup_{t\in \lbrack 0,1]}\left\vert \mathbb{Y}_{n}(\varphi ,t)\right\vert 
\overset{\mathrm{D}}{\longrightarrow }\sup_{t\in \lbrack 0,1]}\left\vert 
\mathbb{Y}(t)\right\vert .
\end{equation*}
\end{theo}

The proof of this theorem is deferred to Section \ref{maintheorem}.
\begin{rem}
\label{rem1} Note that a special case of Theorem \ref{theo1} says that for
each $t\in \left( 0,1\right) $,%
\begin{equation}
\frac{\left( \lfloor nt\rfloor \left( n-\lfloor nt\rfloor \right) \right)
^{2}}{n^{2}(n-1)}U_{\lfloor nt\rfloor ,n}(\varphi )\overset{\mathrm{D}}{%
\longrightarrow }\mathbb{Y}(t).  \label{yy}
\end{equation}%
This fixed $t$ result can be derived from part (a) of Theorem 1.1 of \cite%
{neuhaus1977functional}. \cite{szekely2004testing} point out that
convergence in distribution of a statistic asymptotically equivalent to the
left side of (\ref{yy}) to a nondegenerate random variable should follow
from \cite{neuhaus1977functional} under the null hypothesis of equal
distributions in the two sample case that they consider. Also see \cite%
{rizzo2002test}. (\cite{rizzo2002test} also discuss the consistency of their statistic.) To
the best of our knowledge, we are the first to identify the limit
distribution of the $U_{\lfloor nt\rfloor ,n}(\varphi )$. We should point
out here that the weak convergence result in\ Theorem \ref{theo1} does not
follow from Neuhaus' theorem \cite%
{neuhaus1977functional}, since his result is based on two independent samples,
whereas ours concerns one sample. 
\end{rem}

As suggested in \cite{matteson2014nonparametric}, under the following
assumption, a convergence with probability 1 result can be proved for the
empirical statistic $\mathcal{E}_{k,n}({\bf X};\beta )$ in (\ref{empdiv}). We
shall show that this is indeed the case. 
\begin{ass} 
\label{A}
Let $Y_{i}$, $i\geq1$, and $Z_{i}$, $i\geq1$, be independent i.i.d. sequences, respectively $F_Y$ and $F_Z$.
Also let $Y,Y^{\prime}$ be i.i.d. $F_{Y}$ and $Z,Z^{\prime}$ be i.i.d. $F_{Z}%
$, with $Y,Y^{\prime},Z$ and $Z^{\prime}$ mutually independent. Assume that
for some $\beta\in(0,2)$, $\mathbb{E}(\vert Y\vert ^{\beta
}+\vert Z\vert ^{\beta})<\infty$. Choose $\gamma\in(0,1)$. For
any given $n>1/\gamma$, let $X_{i}=Y_{i}$, for $i=1,\dots,\lfloor
n\gamma\rfloor$, and $X_{i+\lfloor n\gamma\rfloor }=Z_{i}$,
for $i=1,\dots,n-\lfloor n\gamma\rfloor$.
\end{ass}
\begin{lem}
\label{consist} Whenever for a given $\beta \in (0,2)$ Assumption \ref{A}
holds, with probability 1 we have: 
\begin{equation}
\mathcal{E}_{\left\lfloor n\gamma \right\rfloor,n }({\bf X};\beta )\rightarrow 
\mathcal{E}(Y,Z;\beta ).  \label{CON}
\end{equation}
\end{lem}
The proof of this can be found in the Appendix. Next, let $\varphi
(x,y)=|x-y|^{\beta }$, $\beta \in (0,2)$. We see that for any $\gamma \in
(0,1)$ for all large enough $n$, 
\begin{equation*}
\sup_{t\in \lbrack 0,1]}\left\vert \mathbb{Y}_{n}(\varphi ,t)\right\vert
\geq \frac{\left( \lfloor n\gamma \rfloor \left( n-\lfloor n\gamma \rfloor
\right) \right) ^{2}}{n^{2}(n-1)}\mathcal E_{\lfloor n\gamma \rfloor ,n}({\bf X};\beta ),
\end{equation*}%
where it is understood that Assumption \ref{A} holds. Thus by Lemma \ref{consist}, under Assumption \ref{A}, whenever $F_{Y}\neq
F_{Z}$, with probability 1,%
\begin{equation*}
\sup_{t\in \lbrack 0,1]}\left\vert \mathbb{Y}_{n}(\varphi ,t)\right\vert
\rightarrow \infty .
\end{equation*}%
This shows that change-point tests based on the statistic $\sup_{t\in
\lbrack 0,1]}\left\vert \mathbb{Y}_{n}(\varphi ,t)\right\vert $, under the
\textit{sequence of alternatives} of the type given by Assumption \ref{A}, are
consistent. This also has great practical use when looking for
change-points. Intuitively, the $k\in \{1,\hdots,n\}$ that maximizes (\ref%
{empdiv}) would be a good candidate for a change-point location.

\section{From theory to practice\label{practical}}

Theorem \ref{theo1} and the consistency result that follows it lay a firm theoretical foundation to justify the change-point method introduced in \cite%
{matteson2014nonparametric}.
For the present article, since we are not aware of a closed form
expression for the distribution function of the limit process, we may
imagine that this asymptotic result is of limited practical use. Remarkably,
it turns out that we can efficiently approximate via simulation the
distribution of its supremum, leading to a new change-point detection
algorithm with similar performance to \cite{matteson2014nonparametric} but
much faster for longer signals. For instance, finding and testing one change-point in a signal of length $5\,000$ takes eight seconds with
our method and eight minutes using \cite{matteson2014nonparametric}.

To simulate the process $\mathbb{Y}$ we need true or estimated values of the 
$\lambda _{i}$. Recall that these are the eigenvalues of the operator $A$
defined in (\ref{op}). Following \cite{koltchinskii2000random}, the (usually
infinite) spectrum of $A$ can be consistently approximated by the (finite)
spectrum of the empirical $n\times n$ matrix $\tilde{H}_{n}$ whose $(i,j)$%
-th entry is given by 
\begin{equation*}
\tilde{H}_{n}(X_{i},X_{j})=\frac{1}{n}\left( \varphi (X_{i},X_{j})-\mu
(i)-\mu (j)+\eta \right) ,
\end{equation*}%
where $\mu $ is the vector of row means (excluding the diagonal entry) of
matrix $\varphi (X_{i},X_{j})$ and $\eta $ the mean of its upper-diagonal
elements.

In our experience, the $\lambda_{i}$ estimated in this way tend to be quite
accurate for even small $n$. We assert this because upon simulating longer
and longer i.i.d. signals, rapid convergence of the $\lambda_{i}$ is clear.
Furthermore, as there is an exponential drop-off in their magnitude, working
with only a small number (say 20 or 50) of the largest ones appears to be
sufficient for obtaining good results. We illustrate these claims in Section
\ref{applications}. Let us now present our basic algorithm for detecting and
testing for one potential change-point.\newline

\noindent \textbf{Algorithm for detecting and testing one change-point}

\begin{enumerate}
\item Given signal $X_{1},\hdots,X_{n}$, $n\geq 4$, find the $2\leq k\leq
n-2 $ that maximizes the original empirical divergence given in (\ref{empdiv}%
) multiplied by the correct normalization given in (\ref{correctnor}), i.e., 
$k^{2}(n-k)^{2}/n^{2}(n-1)$, and denote the value of this maximum $t^{\star
} $.

\item Calculate the $m$ largest (in absolute value) eigenvalues of the
matrix $\tilde{H}_n$, where $\varphi(X_i,X_j) = |X_i -X_j|^\beta$ and $\beta
\in (0,2)$.

\item Simulate $R$ times the $m$-truncated version of $\mathbb{Y}(t)$ using
the $m$ eigenvalues from the previous step. Record the $R$ values $s_{1},%
\hdots,s_{R}$ of the (absolute) supremum of the process obtained.

\item Reject the null hypothesis of no distributional change (at level $%
\alpha $) if $t_{\mbox{\footnotesize crit}}\leq \alpha $, where 
\hbox{$t_{\mbox{\footnotesize crit}} := \frac{1}{R}\sum_{r=1}^R
\mathbf{1}_{\{s_r > t^{\star}\}}$.} In this case, we deduce a change-point
at the $k$ at which $t^{\star }$ is found. Typically, we set $\alpha =0.05$.
\end{enumerate}
\begin{rem}
One may imagine extending this approach to the multiple change-point case by simply iterating the above algorithm to the left and right of the first-found change-point, and so on. However, as soon as we suppose there can be more than one change-point, the assumption that we may have
 $X_1,\ldots,X_k$ i.i.d., with a different distribution to $X_{k+1},\ldots,X_n$ i.i.d., is immediately broken. Therefore the theory we have presented does not directly follow over to the multiple change-point case. It would be interesting to cleanly extend the results to this, but this would require further theory and multiple testing developments, which are out of the scope of the present article (for references in this direction, see, e.g., \cite{korkas2014multiple}). 
 \end{rem}
The E-divisive algorithm described in \cite{matteson2014nonparametric}
follows a similar logic to our approach except that $t_{\emph{\footnotesize crit}}$ is
calculated via permutation (see \cite{romano2005exact}). Instead of steps 2
and 3, the order of the $n$ data is permuted $R$ times and for the $r$-th
permuted signal, $1 \leq r \leq R$, step 1 is performed to obtain the
absolute maximum $s_r$. The same step 4 is then used to accept or reject the
change-point.

The permutation approach (E-divisive) of \cite{matteson2014nonparametric} is
effective for short signals. Indeed, \cite{hoeffding1952large} showed that
if one can perform all possible permutations, the method produces a test
that is level $\alpha$. However, a signal with $n = 10$ points already
implies more than three million permutations, so a Monte Carlo strategy
(i.e., subsampling permutations with replacement) becomes necessary,
typically with $R=499$. This also gives a test that is theoretically level $%
\alpha$ (see \cite{romano2005exact}) but with much-diminished power.

One could propose increasing the value of $R$ but there is an unfortunate
computational bottleneck in the approach. Usually, one stores in memory the
matrix of $|X_i -X_j|^\beta$ in order to efficiently permute rows/columns
and therefore recalculate $t^{\star}$ each time. But for more than a few
thousand points, manipulating this matrix is slow if not impossible due to
memory constraints. The only alternative to storing and permuting this
matrix is simply to recalculate it each time for each permutation, but this
is very computationally expensive as $n$ increases. Consequently, the
E-divisive approach is only useful for signals up to a few thousand points.

In contrast to this, our algorithm, based on an asymptotic result, risks
underperforming on extremely short signals, and its performance will also
depend on our ability to estimate well the set of largest $\lambda_i$. In
reality though, it works quite well, even on short signals. The matrix with
entries $|X_i -X_j|^\beta$ needs only to be stored once in memory, and all
standard mathematical software (such as Matlab and R) have efficient
functions for finding its largest $m$ eigenvalues (the \texttt{eigs}
function in Matlab and the \texttt{eigs} function in the R package \texttt{%
rARPACK}). Each iteration of the algorithm's simulation step requires
summing the columns of an $m \times T$ matrix of standard normal variables,
where $m$ is the number of $\lambda_i$ retained and $T$ the number of grid
points over which we approximate the Brownian bridge processes between 0 and 1. For $m=50$ and $T=1\,000$ it takes about one
second to perform this $R=499$ times, and is independent of the number of
points in the signal. In contrast, the E-divisive method takes about ten
seconds for $n=1\,000$, one minute for $n=2\,000$, eight minutes for $%
n=5\,000$, etc. One clearly sees the advantage of our approach for longer
signals.

\section{Experimental validation and analysis\label{applications}}

\subsection{Simulated  examples}

It is very important to start with the simplest possible case in order to
demonstrate the fundamental validity of the new method. A basis for
comparison is the E-divisive method from \cite{matteson2014nonparametric}.
Here, we consider signals of length $n \in \{10, 100, 1\,000, 10\,000 \}$
for which either the whole signal is i.i.d. $\mathcal{N}(0,1)$ or else there
is a change-point of height $c \in \{0.1, 0.2, 0.5, 1, 2, 5\}$ after the $%
(n/2)$-th point, i.e., the second half of the signal is i.i.d. $\mathcal{N}%
(c,1)$.

In the former case, we look at the behavior of the Type I error, i.e., the
probability of detecting a change-point when there was none. We have fixed $%
\alpha =0.05$ and want to see how close each method is to this as $n$
increases. In the latter case, we look at the power of the test associated
to each method, i.e., the probability that an actual change-point is
correctly detected as $n$ and $c$ increase. We averaged over $1\,000$
trials. In the following, unless otherwise mentioned we fix $\beta =1$. For
the asymptotic method, the Brownian bridge processes were simulated 499 times;
similarly, for E-divisive we permuted 499 times. Both null distributions
were therefore estimated using the same number of repeats. Note that we did
not test the E-divisive method for $n=10\,000$ because each of the
1\thinspace 000 trials would have taken around two hours to run. All times
given in this paper are with respect to a laptop with a 2.13 GHz Intel Core
2 Duo processor with 4Gb of memory. Results are presented in Figure \ref%
{mean1d}.
\begin{figure}[htb!]
\begin{center}
\includegraphics[height = 5.5cm]{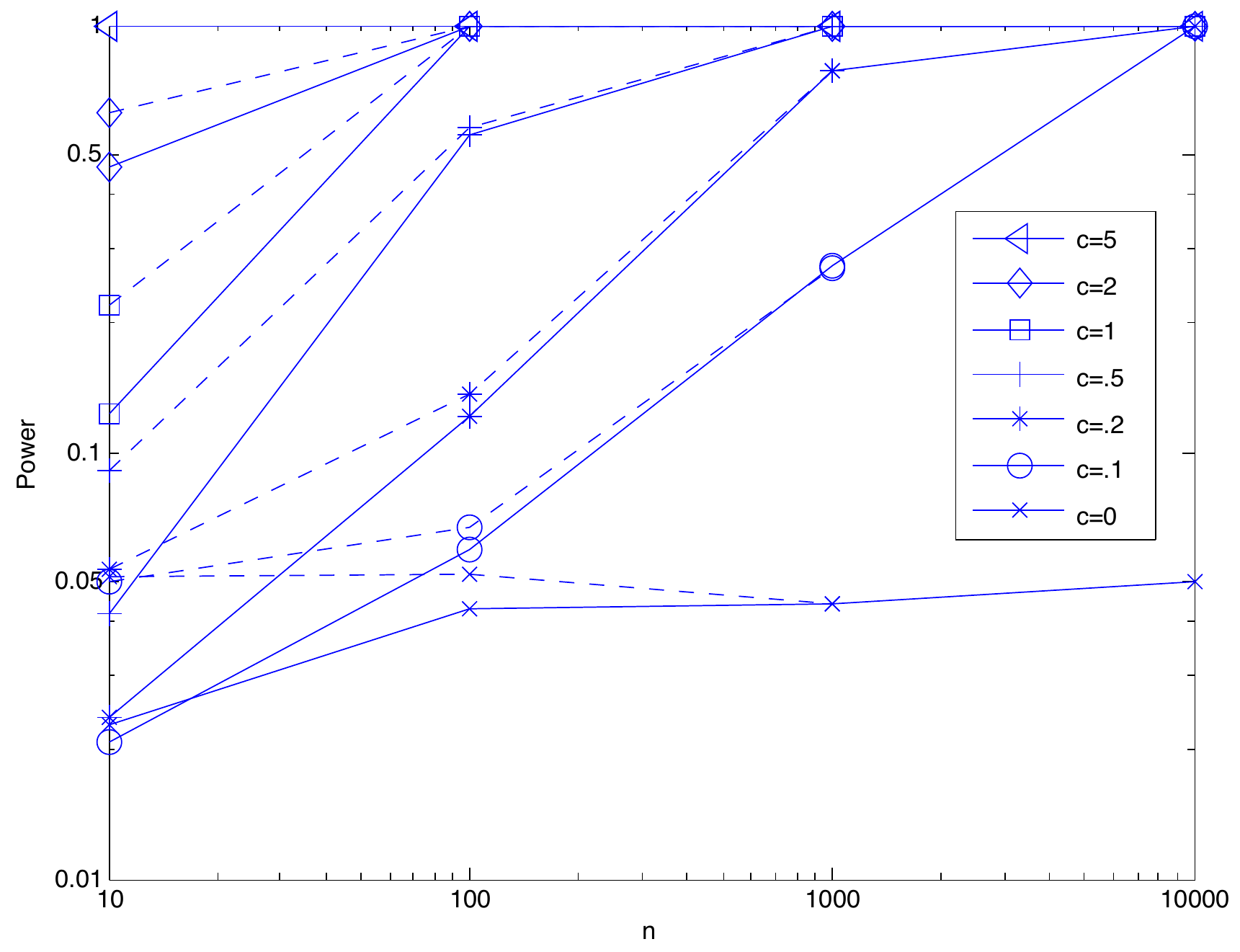}\caption{Statistical power of the asymptotic  (solid line) and E-divisive  (dotted line) methods
for detecting change $c$ in mean  in a Gaussian signal of length $n$. The first $n/2$ points are i.i.d. $\mathcal{N}(0,1)$ and the last $n/2$ points i.i.d. $\mathcal{N}(c,1)$. The Type I error is also shown ($c=0$). Results are averaged over 1\,000 trials.}
\label{mean1d}
\end{center}
\end{figure}

For the Type I error, we see that both methods hover around the intended
value of .05, except for extremely short signals ($n = 10$). As for the
statistical power, it increases as $n$ and $c$ increase. Furthermore, the
asymptotic method rapidly reaches a similar performance as E-divisive: for $%
n=10$, E-divisive is better (but still with quite poor power), for $n=100$
the asymptotic method has almost caught up, and somewhere between $n=100$
and $n=1\,000$ the results become essentially identical; the asymptotic
method has a slight edge at $n=1\,000$.

\begin{figure}[htb!]
\begin{center}
\includegraphics[height = 5.5cm,width = 7.2cm]{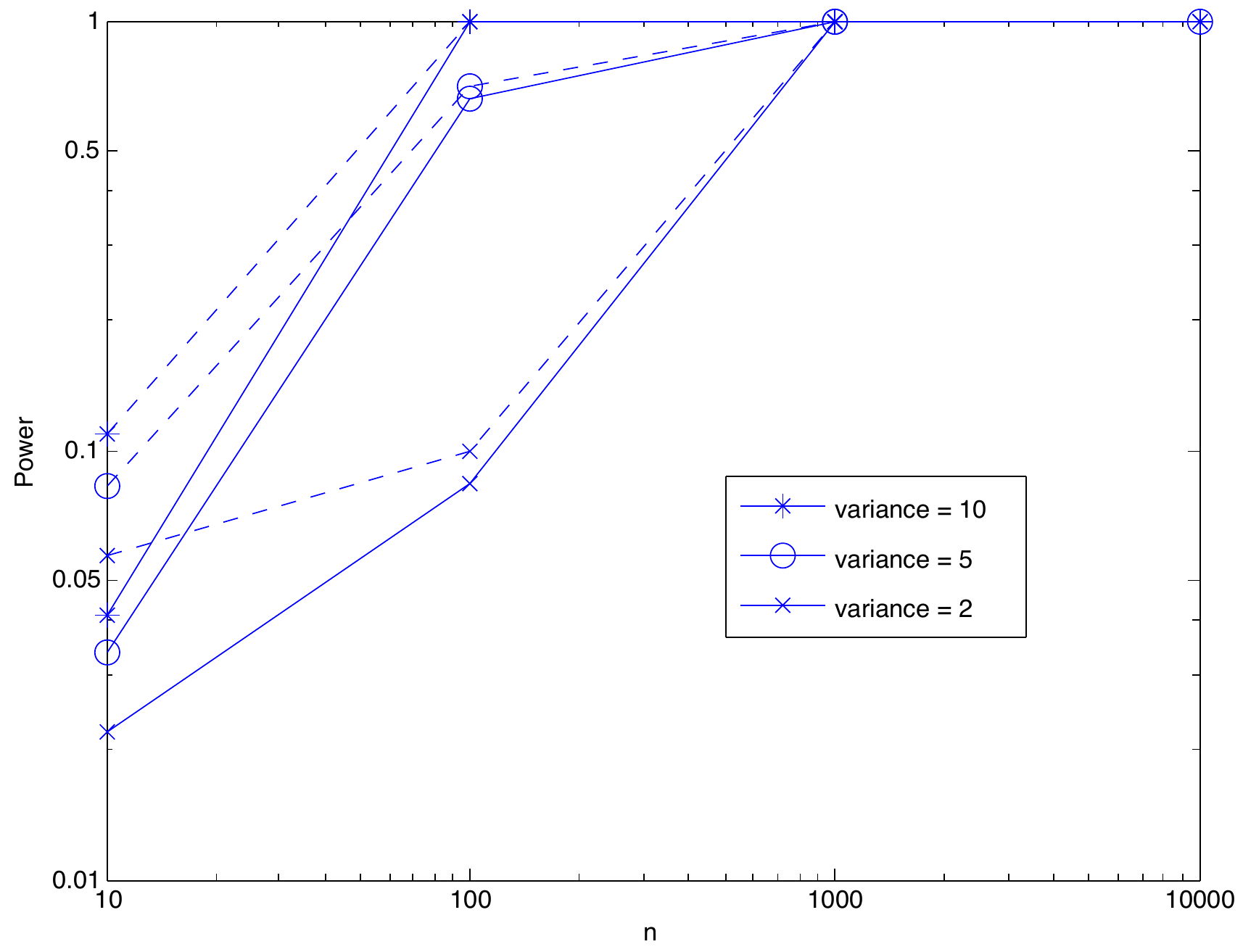}\includegraphics[height = 5.4cm,width = 7.2cm]{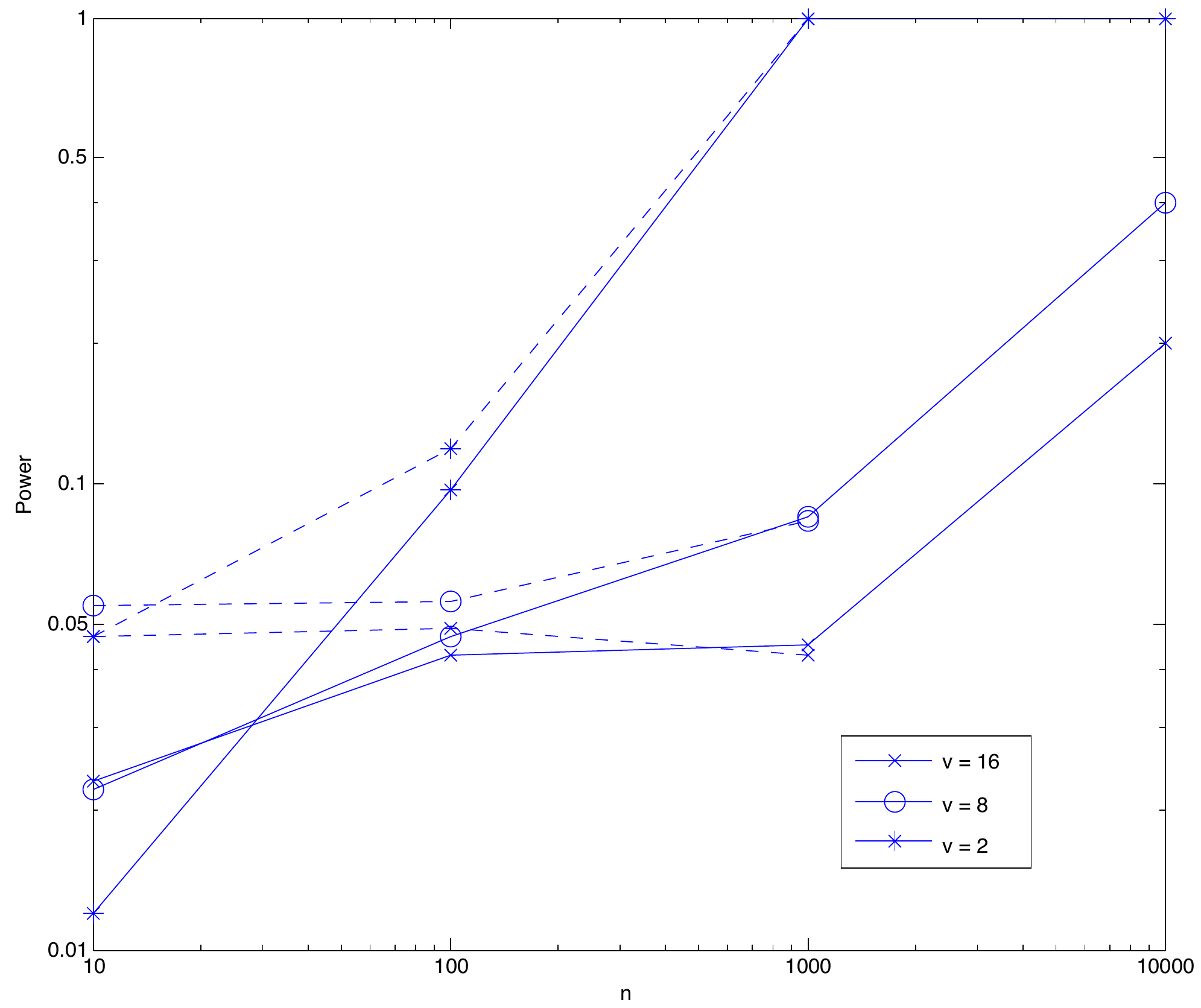}\caption{Statistical power of the asymptotic method (solid line) and E-divisive method (dotted line)
for detecting change in variance (left) and tail (right) in a signal of length $n$. The first $n/2$ points are i.i.d. $\mathcal{N}(0,1)$ and the last $n/2$ points either i.i.d. $\mathcal{N}(0,\sigma^2)$, $\sigma^2 \in \{2, 5, 10\}$ (left) or from a Student's $t_v$ distribution with $v$ degrees of freedom, $v \in \{2, 8, 16\}$ (right). Results are averaged over 1\,000 trials.}
\label{other1d}
\end{center}
\end{figure}

Let us now see to what extent our method is able to detect changes in
variance and tail shape. We considered Gaussian signals of length $n \in
\{10, 100, 1\,000, 10\,000 \}$ for which there is a change-point after the $%
(n/2)$-th point, i.e., the first half of the signal is i.i.d. $\mathcal{N}%
(0,1)$ and the second half either i.i.d. $\mathcal{N}(0,\sigma^2)$ for $%
\sigma^2 \in \{2, 5, 10 \}$ or i.i.d. Student's $t_v$ distributions with $v
\in \{2, 8, 16\}$. Results were averaged over 1\,000 trials and are shown in
Figure \ref{other1d}. 

As before, the statistical power tends to increase as $%
n$ increases and either $\sigma^2$ increases or $v$ decreases. The
asymptotic method matches or beats the performance of E-divisive starting
somewhere between $n=100$ and $n=1\,000$.

\begin{figure}[htb!]
\begin{center}
\includegraphics[width = 12cm]{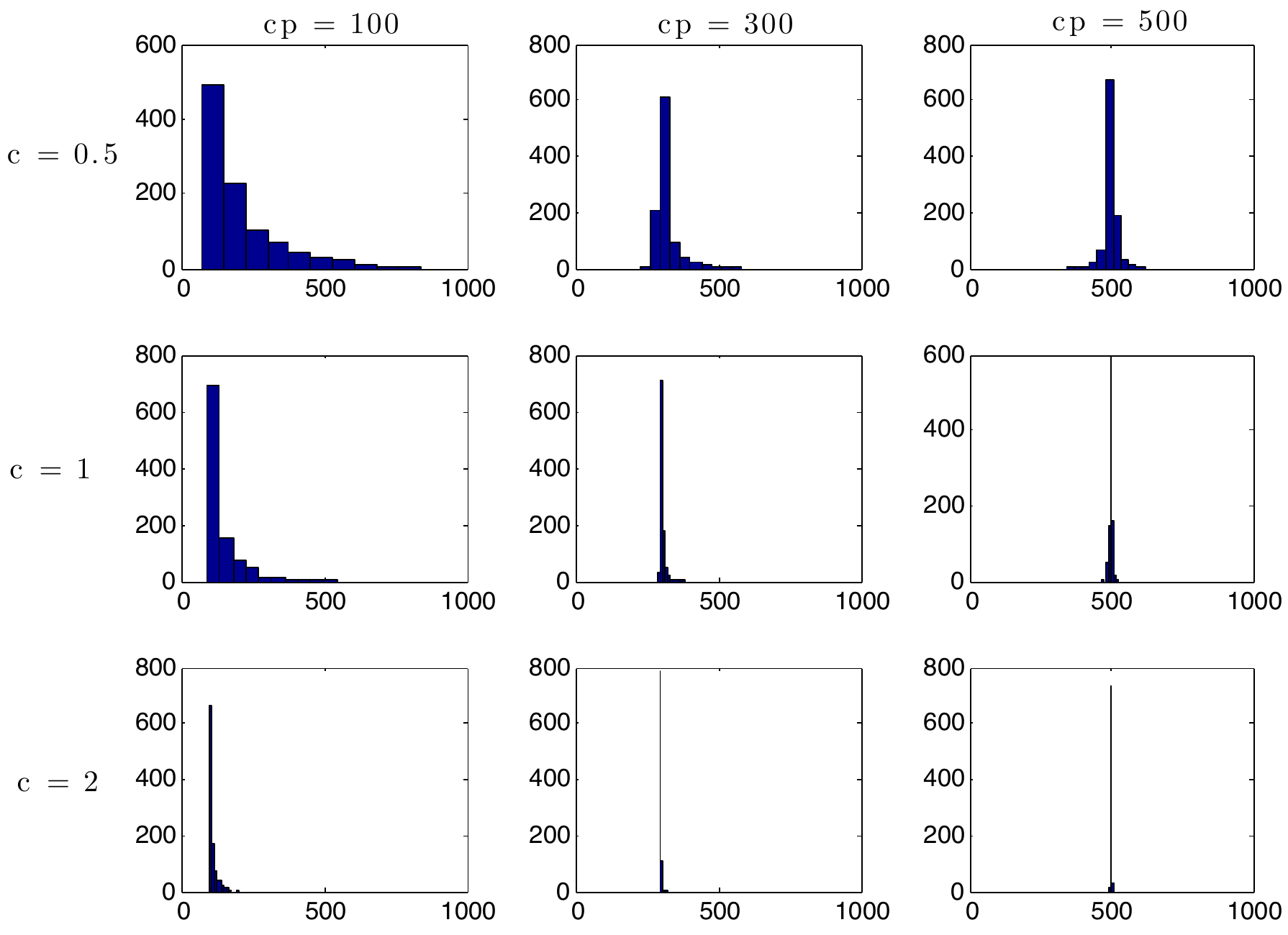}
\caption{Detecting change in mean of $c= 0.5, 1$ or $2$ located at different distances to the boundary (change-point location $\mbox{cp} = 100, 300, 500$) in standardized Gaussian signals with 1\,000 points. Plots show histograms of predicted change-point location over 1\,000 trials.}
\label{move_cp}
\end{center}
\end{figure}

Next, we take a look at the performance of the algorithm when the change-point location moves closer to the boundary. As an illustrative example, we work with sequences of length 1\,000 and either place the change-point after the 100th, 300th or 500th point. Figure \ref{move_cp} shows histograms of 1\,000 repetitions for the predicted location of the change-point, here a change in mean of $c=0.5$ (hardest), $c=1$ (medium) and $c=2$ (easiest). \
We see that moving towards the boundary increases the variance and bias in the prediction. However, as the problem becomes easier (bigger jump in mean), both the variance and bias decrease. Similar results are found when looking at change in variance and tail distribution.
\subsection{Algorithm for long signals}

Remember that as it currently stands, the longest signal that we can treat
depends on the largest matrix that can be stored, which depends in turn on
the memory of a given computer (memory problems for simply manipulating a matrix
on a standard PC typically start to occur around $n=10\,$-$15\,000$). For
this reason, we now propose a modified algorithm that can treat vastly
longer signals. \newline  

\noindent
\textbf{Long-signal algorithm}

\begin{enumerate}

\item Extract sub-signal of equidistant points of length 2\,000. 

\item Run the one change-point algorithm on this. If the null hypothesis is rejected, output the index $k$ of the predicted change-point in this sub-signal. Otherwise, state that no change-point was found.

\item If a change-point was indeed predicted, get the location $k'$ in the original signal corresponding to $k$ in the sub-signal and repeat step 1 of the one change-point algorithm in the interval $[k' - z,\ k' + z]$ to refine the prediction, where $z$ is user-chosen. If $\ell$ is the length of the interval between sub-signal points, one possibility is $z := \min(2\ell, 1\,000)$, where the 1\,000 simply ensures this refining step receives a computationally feasible signal length of at most 2\,000 points.
\end{enumerate}

\begin{figure}[htbp!]
\begin{center}
\includegraphics[height = 8cm]{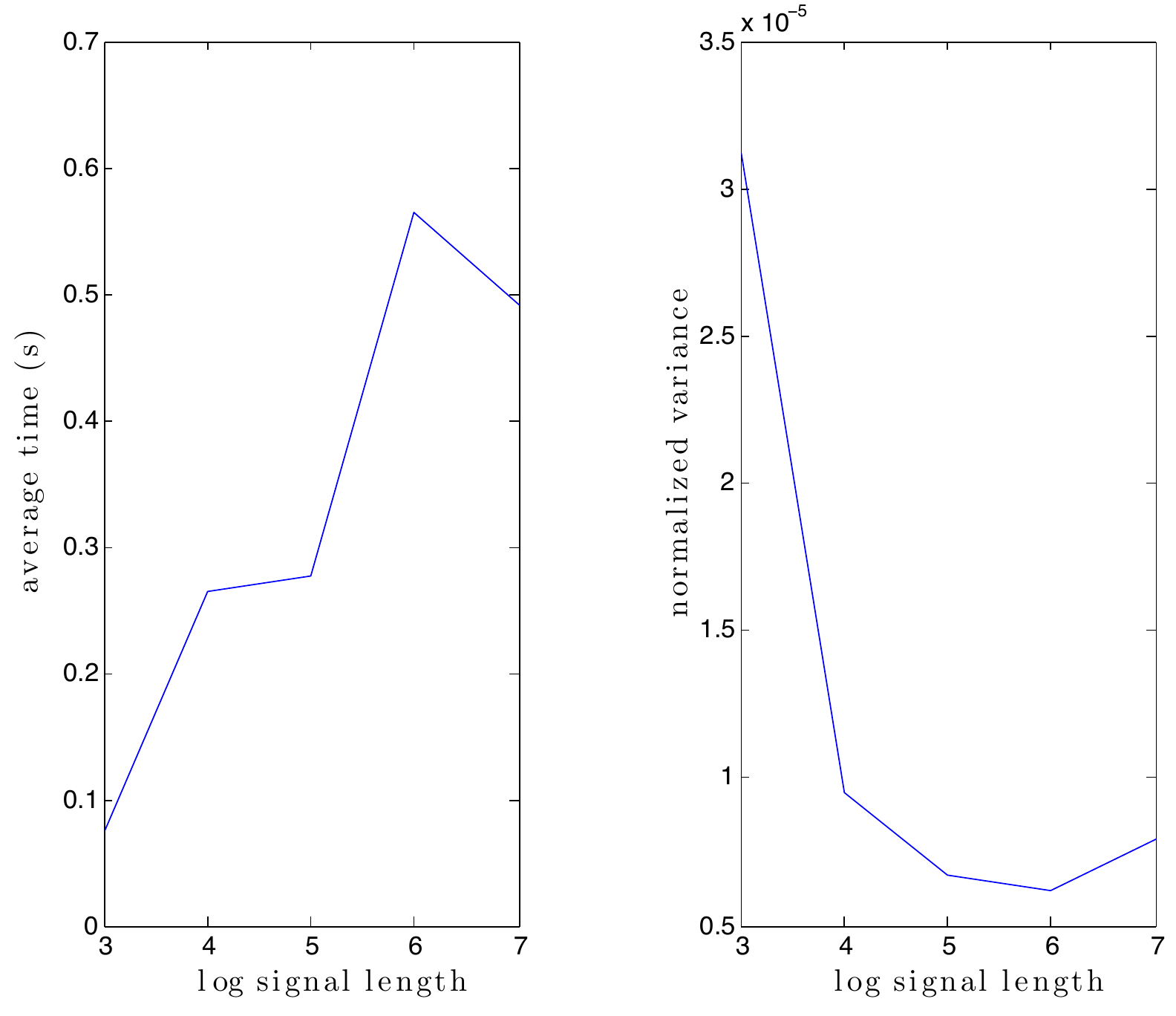}\caption{Long-signal change-point detection. Left: Computing time for signals with 1\,000 to 10 million points. Right: Variance in first change-point prediction over 1\,000 trials after scaling signals to the interval $[0, 1]$.}
\label{big}
\end{center}
\end{figure}

We tested this strategy on simulated standard Gaussian signals of length $10^3,
10^4, 10^5, 10^6$ and $10^7$ with one change-point at the midpoint, a jump of $1$ in the mean. 
Figure \ref{big} (left) shows the time required to locate the potential change-point. 

Clearly, this is  rapid for even extremely long signals. Looking at the algorithm, we see that it merely involves finding a change-point twice, once in the sub-signal, then once in a contiguous block of the original signal of at most length 2\,000. As these two tasks are extremely rapid, the increase in computation time seen mostly comes from the computing overhead of having to extract the sub-signal from longer and longer vectors in memory. 
In Figure \ref{big} (right), we plot the log signal length against the normalized variance, which means that we calculate the variance in predicted change-point location over 1\,000 trials after first dividing the predictions by the length of the signal. Thus all transformed predictions are in the interval 
$[0,1]$ before their variance is taken. This shows that relative to the length of the signal, subsampling does not deteriorate the change-point prediction quality. Instead, what deteriorates due to subsampling is the \emph{absolute} prediction quality, i.e., the variance in predicted change-point location does increase as the signal length increases. However,  we cannot get around this without introducing significantly more sophisticated subsampling procedures, beyond the scope of the work here.

\section{Discussion\label{discussion}}

We have derived the asymptotic distribution of a statistic that was
previously used to build algorithms for finding change-points in signals.
Our new result led to a novel way to construct a practical algorithm for general
change-point detection in long  signals, which came from the surprising realization that it was possible to approximately
simulate from this quite complicated asymptotic distribution. Furthermore,
the method appears to have higher power (in the statistical sense) than
previous methods based on permutation tests for signals of a thousand points
or more. We tested the algorithm on several simulated data
sets, as well as a subsampling variant for dealing with extremely long
signals.

An interesting line of future research would be to find ways to segment the
original signal without requiring stocking a matrix in memory with the same
number of rows and columns as there are points in the signal, currently a
bottleneck for our approach and even more so for previous permutation
approaches. Furthermore, the pertinent choice of the power $%
\beta \in (0,2)$ remains an open question. Lastly, theoretically valid and experimentally feasible extensions of this framework to the multiple change-point case could be a fruitful line of future research.  

\section{Proof of Theorem \protect\ref{theo1}\label{maintheorem}}

To prove Theorem \ref{theo1}, we require a useful technical result. Let us
begin with some notation. For each integer $K\geq 1$, let $D^{K}[0,1]$
denote the space of bounded measurable functions defined on $[0,1]$ taking
values in $\mathbb{R}^{K}$ that are right-continuous with left-hand limits.
For each integer $n\geq 1$, let $\mathbb{V}_{n}^{(k)}$, $k\geq 1$, be a
sequence of processes taking values in $D^{1}[0,1]$ such that for some $M>0$%
, uniformly in $k\geq 1$ and $n\geq 1$,%
\begin{equation}
\mathbb{E}\left( \sup_{t\in \lbrack 0,1]}\left\vert \mathbb{V}%
_{n}^{(k)}(t)\right\vert \right) \leq M.  \label{m1}
\end{equation}%
For each integer $K\geq 1$, define the process taking values in $D^{K}[0,1]$
by 
\begin{equation*}
\mathbb{V}_{n,K}=\left( \mathbb{V}_{n}^{(1)},\hdots,\mathbb{V}%
_{n}^{(K)}\right) .
\end{equation*}%
Assume that for each integer $K\geq 1$, $\mathbb{V}_{n,K}$ converges weakly
as $n\rightarrow \infty $ to the $D^{K}[0,1]$--valued process $\mathbb{V}%
_{K} $ defined as 
\begin{equation*}
\mathbb{V}_{K}:=\left( \mathbb{V}^{(1)},\hdots,\mathbb{V}^{(K)}\right) ,
\end{equation*}%
where $\mathbb{V}^{(k)}$, $k\geq 1,$ is a sequence of $D^{1}[0,1]$--valued
processes such that for some $M>0$, uniformly in $k\geq 1$, 
\begin{equation}
\mathbb{E}\left( \sup_{t\in \lbrack 0,1]}\left\vert \mathbb{V}%
^{(k)}(t)\right\vert \right) \leq M.  \label{m2}
\end{equation}%
We shall establish the following useful result.

\begin{pro}
\label{proposition1} With the notation and assumptions introduced above, for
any choice of constants $a_{m}$, $m\geq1$, satisfying $\sum_{m=1}^{\infty}%
\vert a_{m}\vert <\infty$, the sequence of $D^{1}[ 0,1] $--valued processes 
\begin{equation*}
T_{n}:=\sum_{m=1}^{\infty}a_{m}\mathbb{V}_{n}^{( m) }
\end{equation*}
converges weakly in $D^{1}[ 0,1] $ to the $D^{1}[ 0,1] $--valued process%
\begin{equation*}
T:=\sum_{m=1}^{\infty}a_{m}\mathbb{V}^{( m) }.
\end{equation*}
\end{pro}

\noindent \textit{Proof.} \quad Notice that by (\ref{m1}) 
\begin{equation*}
\mathbb{E}\left( \sum_{m=1}^{\infty }\left\vert a_{m}\right\vert \sup_{t\in
\lbrack 0,1]}\left\vert \mathbb{V}_{n}^{(m)}(t)\right\vert \right) \leq
M\sum_{m=1}^{\infty }|a_{m}|<\infty .
\end{equation*}%
From this we get that with probability $1$, for each $n\geq 1$,%
\begin{equation*}
\sum_{m=1}^{\infty }\left\vert a_{m}\right\vert \sup_{t\in \lbrack
0,1]}\left\vert \mathbb{V}_{n}^{(m)}(t)\right\vert <\infty ,
\end{equation*}%
which in turn implies that with probability $1$, for each $n\geq 1$,%
\begin{equation}
\lim_{K\rightarrow \infty }\sup_{t\in \lbrack 0,1]}\left\vert \overline{T}%
_{n}^{(K)}(t)\right\vert =0,  \label{aa}
\end{equation}%
where 
\begin{equation*}
\overline{T}_{n}^{(K)}(t):=\sum_{m=K+1}^{\infty }a_{m}\mathbb{V}%
_{n}^{(m)}(t).
\end{equation*}%
Since for each $n\geq 1$ and $K\geq 1$, $T_{n}^{(K)}\in D^{1}[0,1]$, where $%
T_{n}^{(K)}:=\sum_{m=1}^{K}a_{m}\mathbb{V}_{n}^{(m)}$, by completeness of $%
D^{1}[0,1]$ in the supremum metric (see page 150 of monograph \cite%
{billingsley1968}), we infer that $T_{n}\in D^{1}[0,1]$. In the same way we
get using (\ref{m2}) that%
\begin{equation}
\lim_{K\rightarrow \infty }\sup_{t\in \lbrack 0,1]}\left\vert \overline{T}%
^{(K)}(t)\right\vert =0,  \label{bb}
\end{equation}%
where 
\begin{equation*}
\overline{T}^{(K)}(t):=\sum_{m=K+1}^{\infty }a_{m}\mathbb{V}^{(m)}(t),
\end{equation*}%
and thus that $T\in D^{1}[0,1]$. Also, since by assumption for each integer $%
K\geq 1$, $\mathbb{V}_{n,K}$ converges weakly as $n\rightarrow \infty $ to
the $D^{K}[0,1]$--valued process $\mathbb{V}_{K}$, we get that $T_{n}^{(K)}$
converges weakly in $D^{1}[0,1]$ to $T^{(K)}$, where 
\begin{equation*}
T_{n}^{(K)}:=\sum_{m=1}^{K}a_{m}\mathbb{V}_{n}^{(m)}\quad \text{ and }\quad
T^{(K)}:=\sum_{m=1}^{K}a_{m}\mathbb{V}^{(m)}.
\end{equation*}%
We complete the proof by combining this with (\ref{aa}) and (\ref{bb}), and
then appealing to Theorem 4.2 of \cite{billingsley1968}. $\quad \square $
\medskip

We are now ready to prove Theorem \ref{theo1}. It turns out that it is more
convenient to prove the result for the following version of the process $%
\mathbb{Y}_{n}$, namely 
\begin{equation*}
\tilde{\mathbb{Y}}_{n}(\tilde{h}_{2},t):=\frac{2\lfloor nt\rfloor \left(
n-\lfloor nt\rfloor \right) }{n^{3}}U_{n}(\tilde{h}_{2})-\frac{2\left(
n-\lfloor nt\rfloor \right) }{n^{2}}U_{\lfloor nt\rfloor }^{(1)}(\tilde{h}%
_{2})-\frac{2\lfloor nt\rfloor }{n^{2}}U_{\lfloor nt\rfloor }^{(2)}(\tilde{h}%
_{2}),
\end{equation*}%
which is readily shown to be asymptotically equivalent to $\mathbb{Y}_{n}(%
\tilde{h}_{2},t)$. Following pages 196-197 of \cite{serfling1980}, we see
that 
\begin{equation*}
\frac{2U_{n}(\tilde{h}_{2})}{n}=\sum_{k=1}^{\infty }\lambda _{k}\left[
\left( \sum_{i=1}^{n}\phi _{k}(X_{i})/\sqrt{n}\right) ^{2}-\frac{1}{n}%
\sum_{i=1}^{n}\phi _{k}^{2}(X_{i})\right] =:\sum_{k=1}^{\infty }\lambda
_{k}\Delta _{k,n},
\end{equation*}%
\begin{equation*}
\frac{2U_{\lfloor nt\rfloor ,n}^{(1)}(\tilde{h}_{2})}{n}=\sum_{k=1}^{\infty
}\lambda _{k}\left[ \left( \sum_{i=1}^{\lfloor nt\rfloor }\phi _{k}(X_{i})/%
\sqrt{n}\right) ^{2}-\frac{1}{n}\sum_{i=1}^{\lfloor nt\rfloor }\phi
_{k}^{2}(X_{i})\right] =:\sum_{k=1}^{\infty }\lambda _{k}\Delta
_{k,n}^{(1)}(t),
\end{equation*}%
and%
\begin{equation*}
\frac{2U_{\lfloor nt\rfloor ,n}^{(2)}(\tilde{h}_{2})}{n}=\sum_{k=1}^{\infty
}\lambda _{k}\left[ \left( \sum_{i=1+\lfloor nt\rfloor }^{n}\phi _{k}(X_{i})/%
\sqrt{n}\right) ^{2}-\frac{1}{n}\sum_{i=1+\lfloor nt\rfloor }^{n}\phi
_{k}^{2}(X_{i})\right] =:\sum_{k=1}^{\infty }\lambda _{k}\Delta
_{k,n}^{(2)}(t).
\end{equation*}%
Thus, 
\begin{equation}
\tilde{\mathbb{Y}}_{n}(\tilde{h}_{2},t)=\sum_{k=1}^{\infty }\lambda
_{k}\left( \frac{\lfloor nt\rfloor \left( n-\lfloor nt\rfloor \right) }{n^{2}%
}\Delta _{k,n}-\frac{\left( n-\lfloor nt\rfloor \right) }{n}\Delta
_{k,n}^{(1)}(t)-\frac{\lfloor nt\rfloor }{n}\Delta _{k,n}^{(2)}(t)\right)
=:\sum_{k=1}^{\infty }\lambda _{k}\mathbb{V}_{n}^{(k)}(t).  \label{v1}
\end{equation}%
Let $\{ \mathbb{W}^{(i)}\} _{i\geq 1}$ be a sequence of standard
Wiener processes on $\left[ 0,1\right] $. Write 
\begin{equation*}
\mathbb{Y}(t):=\sum_{k=1}^{\infty }\lambda _{k}\mathbb{V}^{(k)}(t),
\end{equation*}%
where, for $k\geq 1$,%
\begin{align}
\mathbb{V}^{(k)}(t)& =t(1-t)\left( \left( \mathbb{W}^{(k)}(1)\right)
^{2}-1\right) -(1-t)\left( \left( \mathbb{W}^{(k)}(t)\right) ^{2}-t\right) 
\notag \\
& \quad -t\left( \left( \mathbb{W}^{(k)}(1)-\mathbb{W}^{(k)}(t)\right)
^{2}-(1-t)\right)  \notag \\
& =t(1-t)\left( \left( \mathbb{W}^{(k)}(1)\right) ^{2}+1\right) -(1-t)\left( 
\mathbb{W}^{(k)}(t)\right) ^{2}-t\left( \mathbb{W}^{(k)}(1)-\mathbb{W}%
^{(k)}(t)\right) ^{2}.  \label{v2}
\end{align}%
A simple application of Doob's inequality shows that there exists a constant 
$M>0$ such that (\ref{m1}) and (\ref{m2}) hold, for $\mathbb{V}_{n}^{(k)}$
and $\mathbb{V}^{(k)}$ defined as in (\ref{v1}) and (\ref{v2}).

For any integer $K\geq 1$, let $\mathbb{U}_{1}$ be the random vector such
that $\mathbb{U}_{1}^{T}=(\phi _{1}(X_{1}),\hdots,\phi _{K}(X_{1}))$. We see that
$\mathbb{E}(  \mathbb{U}_{1})  =\mathbf{0}$ and $\mathbb{E}%
(\mathbb{U}_{1}\mathbb{U}_{1}^{T})=I_{K}$.  For any $n\geq 1$ let $\mathbb{U}_{1},\hdots,\mathbb{U}_{n}$ be i.i.d. $%
\mathbb{U}_{1}$. Consider the process defined on $D^{K}[0,1]$
by 
\begin{equation*}
\mathbb{W}_{n,K}(t):=\left( n^{-1/2}\sum_{j\leq \lfloor nt\rfloor }\phi
_{1}(X_{j}),\hdots,n^{-1/2}\sum_{j\leq \lfloor nt\rfloor }\phi
_{K}(X_{j})\right) =:\left( \mathbb{W}_{n}^{(1)}(t),\hdots,\mathbb{W}%
_{n}^{(K)}(t)\right) ,
\end{equation*}%
where for any $i\geq 1$,%
\begin{equation*}
\mathbb{W}_{n}^{(i)}(t):=n^{-1/2}\sum_{j\leq \lfloor nt\rfloor }\phi
_{i}(X_{j}).
\end{equation*}%
Notice that as processes in $t\in [0,1]$,
\[
\mathbb{W}_{n,K}(t):\overset{\mathrm{D}}{=}n^{-1/2}\sum_{j\leq\lfloor
nt\rfloor}\mathbb{U}_{j}.
\]
Clearly by Donsker's theorem the process $(\mathbb{W}_{n,K}(t))_{0\leq t \leq 1}$ converges weakly as $%
n\rightarrow \infty $ to the $\mathbb{R}^{K}$--valued Wiener process $(%
\mathbb{W}_{K}(t))_{0\leq t\leq 1},$ with mean vector zero and covariance
matrix $(t_{1}\wedge t_{2})I_{K}$, $t_{1},t_{2}\in \lbrack 0,1]$, where 
\begin{equation*}
\mathbb{W}_{K}(t):=\left( \mathbb{W}^{(1)}(t),\hdots,\mathbb{W}%
^{(K)}(t)\right) .
\end{equation*}%
Using this fact along with the law of large numbers one readily verifies
that for each integer $K\geq 1$, $(\mathbb{V}_{n}^{(1)},\hdots,\mathbb{V}%
_{n}^{(K)})$ converges weakly as $n\rightarrow \infty $ to $(\mathbb{V}%
^{(1)},\hdots,\mathbb{V}^{(K)})$, where $\mathbb{V}_{n}^{(i)}$ and $\mathbb{V%
}^{(i)}$ are defined as in (\ref{v1}) and (\ref{v2}). All the conditions for
Proposition \ref{proposition1} to hold have been verified. Thus the proof of
Theorem \ref{theo1} is complete, after we note that a little algebra shows
that $\mathbb{Y}(t)$ is equal to 
\begin{equation*}
\sum_{i=1}^{\infty }\lambda _{i}\left( t(1-t) -\left( \mathbb{W}%
^{(i) }( t) -t\mathbb{W}^{(i) }(
1) \right) ^{2}\right) =\sum_{i=1}^{\infty }\lambda _{i}\left( t(
1-t) -\left( \mathbb{B}^{(i) }(t) \right)
^{2}\right) ,  
\end{equation*}%
where $\mathbb{B}^{(i) }(t) =\mathbb{W}^{(
i) }(t) -t\mathbb{W}^{(i) }(1) $, $%
i\geq 1$, are independent Brownian bridges.$\quad \square $

\section{Appendix}

\subsection{Proof of Lemma \protect\ref{consist}}
Notice that for each $n>1$, $\mathcal{E}_{\left\lfloor n\gamma
\right\rfloor,n }({\bf X};\beta )$ is equal to the statistic in (\ref{empdiv}) with $%
k=\left\lfloor n\gamma \right\rfloor $. \ By the law of large numbers for
U-statistics (see Theorem 1 of \cite{sen}) for any $\gamma \in (0,1)$, with
probability 1,
\begin{equation*}
{\binom{\left\lfloor n\gamma \right\rfloor }{2}}^{-1}\sum_{1\leq i<j\leq
\left\lfloor n\gamma \right\rfloor }\left\vert Y_{i}-Y_{j}\right\vert
^{\beta }\rightarrow \mathbb{E}\left\vert Y-Y^{\prime }\right\vert ^{\beta }%
\end{equation*}%
and 
\begin{equation*}
{\binom{n-\left\lfloor n\gamma \right\rfloor }{2}}^{-1}\sum_{1\leq i<j\leq
n-\left\lfloor n\gamma \right\rfloor }\left\vert Z_{i}-Z_{j}\right\vert
^{\beta }\rightarrow \mathbb{E}\left\vert Z-Z^{\prime }\right\vert ^{\beta }.
\end{equation*}%
Next for any $M>0$, write%
\begin{align*}
\left\vert y-z\right\vert ^{\beta } &=\left\vert y-z\right\vert ^{\beta
}1\left\{ \left\vert y\right\vert \leq M,\left\vert z\right\vert \leq
M\right\} +\left\vert y-z\right\vert ^{\beta }1\left\{ \left\vert
y\right\vert \leq M,\left\vert z\right\vert >M\right\}  \\
&\quad +\left\vert y-z\right\vert ^{\beta }1\left\{ \left\vert y\right\vert
>M,\left\vert z\right\vert \leq M\right\} +\left\vert y-z\right\vert ^{\beta
}1\left\{ \left\vert y\right\vert >M,\left\vert z\right\vert >M\right\} .
\end{align*}%
Applying the strong law of large numbers for generalized U-statistics given
in Theorem 1 of \cite{sen}, we get for any $M>0$, with probability 1, 
\begin{align*}
&\frac{2}{\left\lfloor n\gamma \right\rfloor (n-\left\lfloor n\gamma
\right\rfloor )}\sum_{i=1}^{\left\lfloor n\gamma \right\rfloor
}\sum_{j=1}^{n-\left\lfloor n\gamma \right\rfloor }\left\vert
Y_{i}-Z_{j}\right\vert ^{\beta }1\left\{ \left\vert Y_{i}\right\vert \leq
M,\left\vert Z_{j}\right\vert \leq M\right\} \\
&\quad \rightarrow 2\mathbb{E}\left( \left\vert Y-Z\right\vert ^{\beta }1\left\{
\left\vert Y\right\vert \leq M,\left\vert Z\right\vert \leq M\right\}
\right) .
\end{align*}%
Also observe that%
\begin{align*}
&\frac{2}{\left\lfloor n\gamma \right\rfloor (n-\left\lfloor n\gamma
\right\rfloor )}\sum_{i=1}^{\left\lfloor n\gamma \right\rfloor
}\sum_{j=1}^{n-\left\lfloor n\gamma \right\rfloor }\left\vert
Y_{i}-Z_{j}\right\vert ^{\beta }1\left\{ \left\vert Y_{i}\right\vert \leq
M,\left\vert Z_{j}\right\vert >M\right\} \\
&\quad \leq \frac{2}{\left\lfloor n\gamma \right\rfloor (n-\left\lfloor n\gamma
\right\rfloor )}\sum_{i=1}^{\left\lfloor n\gamma \right\rfloor
}\sum_{j=1}^{n-\left\lfloor n\gamma \right\rfloor }\left( M+\left\vert
Z_{j}\right\vert \right) ^{\beta }1\left\{ \left\vert Z_{j}\right\vert
>M\right\}  \\
&\quad=\frac{2}{n-\left\lfloor n\gamma \right\rfloor }\sum_{j=1}^{n-\left\lfloor
n\gamma \right\rfloor }\left( M+\left\vert Z_{j}\right\vert \right) ^{\beta
}1\left\{ \left\vert Z_{j}\right\vert >M\right\} .
\end{align*}%
By the usual law of large numbers for each $M>0$, with probability 1,%
\begin{align*}
\frac{2}{n-\left\lfloor n\gamma \right\rfloor }\sum_{j=1}^{n-\left\lfloor
n\gamma \right\rfloor }\left( M+\left\vert Z_{j}\right\vert \right) ^{\beta
}1\left\{ \left\vert Z_{j}\right\vert >M\right\}  &\rightarrow 2\mathbb{E}%
\left( \left( M+\left\vert Z\right\vert \right) ^{\beta }1\left\{ \left\vert
Z\right\vert >M\right\} \right)  \\
&\leq 2^{\beta +1}\mathbb{E}\left( \left\vert Z\right\vert ^{\beta
}1\left\{ \left\vert Z\right\vert >M\right\} \right) .
\end{align*}%
Thus, with probability 1, for all $M>0$,
\begin{align*}
&\limsup_{n\rightarrow \infty }\frac{2}{\left\lfloor n\gamma \right\rfloor
(n-\left\lfloor n\gamma \right\rfloor )}\sum_{i=1}^{\left\lfloor n\gamma
\right\rfloor }\sum_{j=1}^{n-\left\lfloor n\gamma \right\rfloor }\left\vert
Y_{i}-Z_{j}\right\vert ^{\beta }1\left\{ \left\vert Y_{i}\right\vert \leq
M,\left\vert Z_{j}\right\vert >M\right\}  \\
&\quad \leq 2^{\beta +1}\mathbb{E}\left( \left\vert Z\right\vert ^{\beta
}1\left\{ \left\vert Z\right\vert >M\right\} \right) .
\end{align*}%
In the same way we get that, with probability 1,%
\begin{align*}
&\limsup_{n\rightarrow \infty }\frac{2}{\left\lfloor n\gamma \right\rfloor
(n-\left\lfloor n\gamma \right\rfloor )}\sum_{i=1}^{\left\lfloor n\gamma
\right\rfloor }\sum_{j=1}^{n-\left\lfloor n\gamma \right\rfloor }\left\vert
Y_{i}-Z_{j}\right\vert ^{\beta }1\left\{ \left\vert Y_{i}\right\vert
>M,\left\vert Z_{j}\right\vert \leq M\right\}  \\
&\quad \leq 2\mathbb{E}\left( \left( \left\vert Y\right\vert +M\right) ^{\beta
}1\left\{ \left\vert Y\right\vert >M\right\} \right) \leq 2^{\beta +1}%
\mathbb{E}\left( \left\vert Y\right\vert ^{\beta }1\left\{ \left\vert
Y\right\vert >M\right\} \right) .
\end{align*}%
Finally, note that, by the $c_{r}$-inequality,
\begin{align*}
&\frac{2}{\left\lfloor n\gamma \right\rfloor (n-\left\lfloor n\gamma
\right\rfloor )}\sum_{i=1}^{\left\lfloor n\gamma \right\rfloor
}\sum_{j=1}^{n-\left\lfloor n\gamma \right\rfloor }\left\vert
Y_{i}-Z_{j}\right\vert ^{\beta }1\left\{ \left\vert Y_{i}\right\vert
>M,\left\vert Z_{j}\right\vert >M\right\}\\
&\quad \leq \frac{2^{\beta }}{\left\lfloor n\gamma \right\rfloor (n-\left\lfloor
n\gamma \right\rfloor )}\sum_{i=1}^{\left\lfloor n\gamma \right\rfloor
}\sum_{j=1}^{n-\left\lfloor n\gamma \right\rfloor }\left( \left\vert
Y_{i}\right\vert ^{\beta }+\left\vert Z_{j}\right\vert ^{\beta }\right)
1\left\{ \left\vert Y_{i}\right\vert >M,\left\vert Z_{j}\right\vert
>M\right\} \\
&\quad \leq \frac{2^{\beta }}{\left\lfloor n\gamma \right\rfloor }%
\sum_{i=1}^{\left\lfloor n\gamma \right\rfloor }\left\vert Y_{i}\right\vert
^{\beta }1\left\{ \left\vert Y_{i}\right\vert >M\right\} +\frac{2^{\beta }}{%
n-\left\lfloor n\gamma \right\rfloor }\sum_{j=1}^{n-\left\lfloor n\gamma
\right\rfloor }\left\vert Z_{j}\right\vert ^{\beta }1\left\{ \left\vert
Z_{j}\right\vert >M\right\}.
\end{align*}%
By the law of large numbers this converges, with probability 1, to 
\begin{equation*}
2^{\beta }\mathbb{E}\left( \left\vert Y\right\vert ^{\beta }1\left\{
\left\vert Y\right\vert >M\right\} \right) +2^{\beta }\mathbb{E}\left(
\left\vert Z\right\vert ^{\beta }1\left\{ \left\vert Z\right\vert >M\right\}
\right) .
\end{equation*}%
Obviously as $M\rightarrow \infty $, 
\begin{equation*}
2\mathbb{E}\left( \left\vert Y-Z\right\vert ^{\beta }1\left\{ \left\vert
Y\right\vert \leq M,\left\vert Z\right\vert \leq M\right\} \right)
\rightarrow 2\mathbb{E}\left\vert Y-Z\right\vert ^{\beta }
\end{equation*}%
and 
\begin{equation*}
3\cdot 2^{\beta }\mathbb{E}\left( \left\vert Y\right\vert ^{\beta }1\left\{
\left\vert Y\right\vert >M\right\} \right) +3\cdot 2^{\beta }\mathbb{E}%
\left( \left\vert Z\right\vert ^{\beta }1\left\{ \left\vert Z\right\vert
>M\right\} \right) \rightarrow 0.
\end{equation*}%
Putting everything together we get that (\ref{CON}) holds. $\quad \square $
\subsection{A technical result}
Let $X$ and $X^{\prime }$ be i.i.d. $F$ and let $\varphi $ be a symmetric
measurable function from $\mathbb{R}^{d}\times \mathbb{R}^{d}\rightarrow 
\mathbb{R}$ such that $\mathbb{E}\varphi ^{2}(X,X^{\prime })<\infty $.
Recall the notation (\ref{H1}). Let $A$ be the operator defined on $L_{2}(%
\mathbb{R}^{d},F)$ as in (\ref{op}).

Notice that 
\begin{equation*}
\mathbb{E}\big( g(X)\tilde{h}_{2}(X,X^{\prime })g(X^{\prime })\big) =\int_{%
\mathbb{R}^{d}}g(x)Ag(x)dF(x)=:\langle g,Ag\rangle .
\end{equation*}

Let us now introduce some useful definitions. Given $\beta \in (0,2)$ and $%
\varphi _{\beta }(x,y)=|x-y|^{\beta }$, define as in (\ref{H1}), 
\begin{equation*}
h_{2,\beta }(x,y)=\varphi _{\beta }(x,y)-\varphi _{1,\beta }(x)-\varphi
_{1,\beta }(y)\quad \text{ and }\quad \tilde{h}_{2,\beta }(x,y)=h_{\beta
}(x,y)+\mathbb{E}\varphi _{\beta }(X,X^{\prime }).
\end{equation*}%
The aim here is to verify that the function $\tilde{h}_{2,\beta }(x,y)$
satisfies the conditions of Theorem \ref{theo1} as long as 
\begin{equation}
\mathbb{E}|X|^{2\beta }<\infty .  \label{beta}
\end{equation}%
Let $\tilde{A}_{\beta }$ denote the integral operator 
\begin{equation*}
\tilde{A}_{\beta }g(x)=\int_{\mathbb{R}^{d}}\text{ }\tilde{h}_{2,\beta
}(x,y)g(y)dF(y),\quad x\in \mathbb{R}^{d},\text{ }g\in L_{2}(\mathbb{R}%
^{d},F).
\end{equation*}%
Clearly (\ref{beta}) implies (\ref{sec}) with $\varphi =\varphi _{\beta }$,
which, in turn, by (\ref{lam}) implies 
\begin{equation*}
\mathbb{E}\tilde{h}_{2,\beta }^{2}(X,X^{\prime })=\int_{\mathbb{R}^{d}}\int_{%
\mathbb{R}^{d}}\tilde{h}_{2,\beta }^{2}(x,y)dF(x)dF(y)=\sum_{i=1}^{\infty
}\lambda _{i}^{2}<\infty ,  
\end{equation*}%
where $\lambda _{i}$, $i\geq 1$, are the eigenvalues of the operator $\tilde{%
A}_{\beta }$, with corresponding orthonormal eigenfunctions $\phi _{i}$, $%
i\geq 1$. 

Next we shall prove that when (\ref{beta}) holds then the eigenvalues $%
\lambda _{i}$, $i\geq 1$, of this integral operator $\tilde{A}_{\beta }$
satisfy (\ref{finite}). This is summarized in the following lemma, whose
proof is postponed to the next paragraph.

\begin{lem}
\label{alpha} Whenever for some $\beta \in (0,2)$, (\ref{beta}) holds, the
eigenvalues $\lambda _{i}$, $i\geq 1$, of the operator $\tilde{A}_{\beta }$
satisfy (\ref{finite}).
\end{lem}

The technical results that follow will imply that $\lambda _{i}\leq 0$ for
all $i\geq 1$ and $\sum_{i=1}^{\infty }\lambda _{i}$ is finite, from which
we can infer (\ref{finite}), and thus Lemma \ref{alpha}. Let us begin with
two definitions.

\begin{defi}
\label{definition1} Let $\mathcal{X}$ be a nonempty set. A symmetric
function $K:\mathcal{X}\times \mathcal{X}\rightarrow \mathbb{R}$ is called
positive definite if 
\begin{equation*}
\sum_{i=1}^{n}\sum_{j=1}^{n}c_{i}c_{j}K(x_{i},x_{j})\geq 0
\end{equation*}%
for all $n\geq 1$, $c_{1},\hdots,c_{n}\in \mathbb{R}$ and $x_{1},\hdots%
,x_{n}\in \mathcal{X}$.
\end{defi}

\begin{defi}
\label{definition2}Let $\mathcal{X}$ be a nonempty set. A symmetric function 
$K:\mathcal{X}\times \mathcal{X}\rightarrow \mathbb{R}$ is called
conditionally negative definite if 
\begin{equation*}
\sum_{i=1}^{n}\sum_{j=1}^{n}c_{i}c_{j}K(x_{i},x_{j})\leq 0
\end{equation*}%
for all $n\geq 1$, $c_{1},\hdots,c_{n}\in \mathbb{R}$ such that \ $%
\sum_{i=1}^{n}c_{i}=0$ and $x_{1},\hdots,x_{n}\in \mathcal{X}$.
\end{defi}

Next, we shall be using part of Lemma 2.1 on page 74 of \cite{bergharmonic},
which we state here for convenience as Lemma \ref{lemmaA}.

\begin{lem}
\label{lemmaA} Let $K$ be a symmetric function on $\mathcal{X}\times 
\mathcal{X}$. Then, for any $x_{0}\in \mathcal{X}$, the function 
\begin{equation*}
\tilde{K}(x,y)=K(x,x_{0})+K(y,x_{0})-K(x,y)-K(x_{0},x_{0})
\end{equation*}%
is positive definite if and only if $K$ is conditionally negative definite.
\end{lem}

The following lemma can be proved just as Corollary 2.1 in \cite{FerMen}.
\begin{lem}
\label{lemmaB} Let $H:\mathbb{R}^{d}\times \mathbb{R}^{d}$ $\rightarrow 
\mathbb{R}$ be a symmetric positive definite function in the sense of
Definition \ref{definition1}. Assume that $H$ is continuous and $\mathbb{E}
H^{2}(X,X^{\prime })<\infty ,$ where $X$ and $X^{\prime }$ are i.i.d. $F$.
Then $\mathbb{E}(g(X)H(X,X^{\prime })g(X^{\prime }))\geq 0$ for all $g\in
L_{2}(\mathbb{R}^{d},F)$, i.e., $H$ is $L^{2}$-positive definite in the
sense of \cite{FerMen}.
\end{lem}

We recall that an operator $L$ on $L_{2}(\mathbb{R}^{d},F)$ is called
positive definite if for all $g\in L_{2}(\mathbb{R}^{d},F)$, $\langle
g,Lg\rangle \geq 0$.
\begin{pro}
\label{proposition2} Let $\varphi :\mathbb{R}^{d}\times \mathbb{R}%
^{d}\rightarrow \mathbb{R}$ be a symmetric continuous function that is a
conditionally negative definite function in the sense of Definition \ref%
{definition2}. Assume that $\varphi (x,x)=0$ for all $x\in \mathbb{R}^{d}$
and $\mathbb{E}\varphi ^{2}(X,X^{\prime })<\infty $. Then $\varphi $ defines
a positive definite operator $L$ on $L_{2}(\mathbb{R}^{d},F)$ given by 
\begin{equation*}
Lg(x)=-\int_{\mathbb{R}^{d}}h(x,y)g(y)dF(y),\quad x\in \mathbb{R}^{d},g\in
L_{2}(\mathbb{R}^{d},F),
\end{equation*}%
where $h$ is defined as in (\ref{H1}). Furthermore the operator $\tilde{L}$
on $L_{2}(\mathbb{R}^{d},F)$ given by%
\begin{equation*}
\tilde{L}g(x)=-\int_{\mathbb{R}^{d}}\left( h(x,y)+\mathbb{E}\varphi
(X,X^{\prime })\right) g(y)dF(y),\quad x\in \mathbb{R}^{d},g\in L_{2}(%
\mathbb{R}^{d},F),
\end{equation*}%
is also a positive definite operator on $L_{2}(\mathbb{R}^{d},F)$.
\end{pro}

\textit{Proof.} We must show that for all $g\in L_{2}(\mathbb{R}^{d},F)$, 
\begin{equation*}
\langle g,Lg\rangle =-\mathbb{E}(g(X)h(X,X^{\prime })g(X^{\prime }))\geq 0.
\end{equation*}%
For any $u\in \mathbb{R}^{d}$, let us write 
\begin{align*}
\varphi (x,y,u)& :=\varphi (x,u)+\varphi (y,u)-\varphi (u,u)-\varphi (x,y) \\
& =\varphi (x,u)+\varphi (y,u)-\varphi (x,y).
\end{align*}%
Since $\varphi $ is assumed to be conditionally negative definite, by Lemma %
\ref{lemmaA} we have that for any fixed $u\in \mathbb{R}^{d}$, $\varphi
(x,y,u)$ is positive definite in the sense of Definition \ref{definition1}.
Hence, since $\varphi $ is also assumed to be continuous, by Lemma \ref%
{lemmaB} for all $g\in L_{2}(\mathbb{R}^{d},F)$, 
\begin{equation*}
\mathbb{E}\left( g(X)\varphi (X,X^{\prime },u)g(X^{\prime })\right) \geq 0.
\end{equation*}%
Noting that if $U$ has distribution function $F$, $\mathbb{E}\varphi
(x,y,U)=-h(x,y)$, we get, assuming that $X$, $X^{\prime }$ and $U$ are
independent, that 
\begin{equation*}
\mathbb{E}\left( g(X)\varphi (X,X^{\prime },U)g(X^{\prime })\right) =-%
\mathbb{E}\left( g(X)h(X,X^{\prime })g(X^{\prime })\right) \geq 0.\qquad
\end{equation*}%
Next, notice that for any eigenvalue and normalized eigenfunction $(\tilde{%
\lambda}_{i},\tilde{\phi}_{i})$ pair, $i\geq 1$, of the operator $%
\tilde{L}$, we have 
\begin{equation*}
\tilde{\lambda}_{i}\tilde{\phi}_{i}(x)=\tilde{L}\tilde{\phi}_{i}(x)=-\int_{%
\mathbb{R}^{d}}\left( h(x,y)+\mathbb{E}\varphi (X,X^{\prime })\right) \tilde{%
\phi}_{i}\left( y\right) dF(y).
\end{equation*}%
Now, 
\begin{equation*}
\int_{\mathbb{R}^{d}}\left( h(x,y)+\mathbb{E}\varphi (X,X^{\prime })\right)
dF(y)=0\text{, for all }x\in \mathbb{R}^{d},
\end{equation*}%
implies that $(\tilde{\lambda}_{1}, \tilde{\phi}_{1}):=\left(
0,1\right) $ is an eigenvalue and normalized eigenfunction pair of $\tilde{L}$. From this
we get that whenever $\tilde{\lambda}_{i}\neq 0$, $\mathbb{E}\tilde{\phi}%
_{i}(X)=0$, $i\geq 2$, which says that for such $\tilde{\lambda}_{i}$, 
\begin{equation*}
\tilde{\lambda}_{i}\tilde{\phi}_{i}(x)=-\int_{\mathbb{R}^{d}}h(x,y)\tilde{%
\phi}_{i}(y)dF(y).
\end{equation*}%
This implies that whenever for some $i\geq1,$ $( \tilde{\lambda}_{i},\tilde{\phi}_{i})$, with $\tilde{\lambda}_{i}\neq0,$ is an eigenvalue and
normalized eigenfunction pair of the operator $\tilde{L},$ it is also an
eigenvalue and normalized eigenfunction pair of the operator $L$. Moreover,
since the integral operator $L$ is positive definite on $L_{2}(\mathbb{R}%
^{d},F)$, this implies that for any such nonzero $\tilde{\lambda}_{i}$
(where necessarily $i\geq2$)%
\begin{align*}
&-\int_{\mathbb{R}^{d}}\int_{\mathbb{R}^{d}}\tilde{\phi}_{i}(x)\left( h(x,y)+%
\mathbb{E}\varphi(X,X^{\prime})\right) \tilde{\phi}_{i}(y)dF(x)dF(y) \\
&\quad =-\int_{\mathbb{R}^{d}}\int_{\mathbb{R}^{d}}\tilde{\phi}_{i}(x)h(x,y)%
\tilde{\phi}_{i}(y)dF(x)dF(y)=\tilde{\lambda}_{i}\geq0,
\end{align*}
which says that the operator $\tilde{L}$ is positive definite on $L_{2}( 
\mathbb{R}^{d},F)$. $\quad \square$

\subsection{Proof of Lemma \protect\ref{alpha}}

A special case of Theorem 3.2.2 in \cite{bergharmonic} says that the
function $\varphi _{\beta }(x,y)=|x-y|^{\beta }$, $\beta \in (0,2)$, is
conditionally negative definite. Also see Exercise 3.2.13b in \cite%
{bergharmonic} and the discussion after Proposition 3 in \cite%
{szekely2013energy}. Therefore by Proposition \ref{proposition2} the
integral operator $L_{\beta }$ defined by the function 
\begin{equation*}
K_{\beta }\left( x,y\right) =-h_{2,\beta }\left( x,y\right)
\end{equation*}%
is positive definite as well as the integral operator $\tilde{L}_{\beta } =-%
\tilde{A}_{\beta }$ defined by the function 
\begin{equation*}
\tilde{K}_{\beta }\left( x,y\right) =-h_{2,\beta }(x,y)-\mathbb{E}\varphi
_{\beta }(X,X^{\prime }).
\end{equation*}%

Next, as in the proof of Proposition \ref{proposition2}, any eigenvalue and
normalized eigenfunction 
\begin{equation*}
\big( \tilde{\lambda }_{i},\tilde{\phi }_{i}\big) =\left( -\lambda
_{i},-\phi _{i}\right)
\end{equation*}%
pair, with $\tilde{\lambda }_{i}\neq 0$, $i\geq 1$, of the operator $\tilde{L%
}_{\beta }=-\tilde{A}_{\beta }$ is also an eigenvalue and normalized
eigenfunction pair of the operator $L_{\beta }$. 

We shall apply Theorem 2 of \cite{sun} to show that uniformly on compact
subsets $D$ of $\mathbb{R}^{d}$, 
\begin{equation*}
K_{\beta }\left( x,y\right) =\sum_{i=1}^{\infty }\rho _{i}\psi _{i}(x)\psi
_{i}(y),\quad \left( x,y\right) \in D\times D,
\end{equation*}%
where $\rho _{i}\geq 0$, $i\geq 1$, are the eigenvalues of the operator $%
L_{\beta }=-A_{\beta }$ with corresponding normalized eigenfunctions $\psi
_{i}$, $i\geq 1$. In particular
\begin{equation*}
K_{\beta }\left( x,x\right) =\sum_{i=1}^{\infty }\rho _{i}\psi _{i}^{2}(x), \quad x\in D,
\end{equation*}%
and thus since $\mathbb{E}\psi _{i}^{2}(X)=1$, $i\geq 1$, and $\mathbb{E}%
K_{\beta }\left( X,X\right) <\infty $, we get 
\begin{equation*}
\sum_{i=1}^{\infty }\rho _{i}<\infty .
\end{equation*}%
Therefore since, as pointed out above, the eigenvalue and normalized eigenfunction pairs $%
( -\lambda_{i},-\phi_{i}) $ of $\tilde{L}_{\beta }=-\tilde{A}_{\beta}$, with 
$\lambda_{i}\neq0$, are also eigenvalue and normalized eigenfunction pairs of the operator 
$L_{\beta}$ this implies that $\sum _{i=1}^{\infty}\vert \lambda_{i}\vert
<\infty$.

Our proof will be complete once we have checked that $L_{\beta }$ satisfies
the conditions of Theorem 2 of \cite{sun}.

Since the function $\varphi _{\beta }(x,y)=|x-y|^{\beta }$, $\beta \in
(0,2), $ is conditionally negative definite, by Lemma \ref{lemmaA} the
function $K_{\beta }\left( x,y\right) $ is positive definite. To see this
note that by Lemma \ref{lemmaA} for any fixed $u\in\mathbb{R}$ the function%
\begin{equation*}
\varphi_{\beta}\left( x,u\right) +\varphi_{\beta}\left( u,y\right)
-\varphi_{\beta}\left( x,y\right) -\varphi_{\beta}\left( u,u\right)
=\varphi_{\beta}\left( x,u\right) +\varphi_{\beta}\left( u,y\right)
-\varphi_{\beta}\left( x,y\right)
\end{equation*}
is positive definite. Therefore we readily see that 
\begin{equation*}
K_{\beta}(x,y)=\left( \int_{\mathbb{R}^{d}}\varphi_{\beta}\left( x,u\right)
+\varphi_{\beta}\left( u,y\right) -\varphi_{\beta}\left( x,y\right) \right)
dF\left( u\right)
\end{equation*}
is positive definite. In addition, $K_{\beta }\left( x,y\right) $ is
symmetric and continuous, and thus $K_{\beta }\left( x,y\right) $ is a
Mercer kernel in the terminology of \cite{sun}. We must also verify the
following assumptions. $\smallskip $

\noindent \textbf{Assumption A}. For each $x\in \mathbb{R}^{d}$, $K_{\beta
}\left( x,\cdot \right) \in L_{2}(\mathbb{R}^{d},F).$ $\smallskip $

\noindent \textbf{Assumption B}. $L_{\beta }$ is a bounded and positive
definite operator on $L_{2}(\mathbb{R}^{d},F)$ and for every $g\in L_{2}(%
\mathbb{R}^{d},F)$, the function $\ $%
\begin{equation*}
L_{\beta }g( x) =\int_{\mathbb{R}^{d}}K_{\beta }(x,y)g(y)dF(y)
\end{equation*}%
is a continuous function on $\mathbb{R}^{d}$.$\smallskip $

\noindent \textbf{Assumption C}. $L_{\beta }$ has at most countably many
positive eigenvalues and orthonormal eigenfunctions.

Since $\varphi _{\beta }$ is a symmetric continuous function that is
conditionally negative definite in the sense of Definition \ref{definition2}
satisfying $\varphi _{\beta }(x,x)=0$ for all $x\in \mathbb{R}^{d}$ and $%
\mathbb{E}\varphi _{\beta }^{2}(X,X^{\prime })<\infty $, we get by
Proposition \ref{proposition2} that $L_{\beta }$ is a positive definite
operator on $L_{2}(\mathbb{R}^{d},F)$. \ Also (\ref{beta}) obviously implies
that Assumption A holds and 
\begin{equation*}
\mathbb{E}K_{\beta }^{2}(X,X^{\prime })=\int_{\mathbb{R}^{d}}\int_{\mathbb{R}%
^{d}}K_{\beta }^{2}(x,y)dF(x)dF(y)<\infty ,  
\end{equation*}%
which by Proposition 1 of \cite{sun} implies that the operator $L_{\beta }$
is bounded and compact. (From Sun's Proposition 1 one can also infer that $%
L_{\beta }$ is positive definite. However, he does not provide a proof.
Therefore we invoke our Lemma \ref{lemmaB} here.) An elementary argument
based on the dominated convergence theorem implies that $L_{\beta } g ( x) $
is a continuous function on $\mathbb{R}^{d}$. Thus Assumption B is
satisfied. Finally, since the operator $L_{\beta }$ is compact, Theorem
VII.4.5 of \cite{dunford1988linear} implies that Assumption C is fulfilled.
Thus the assumptions of Theorem 2 of \cite{sun} hold. This completes the
proof of Lemma \ref{alpha}. $\quad \square $

\bibliographystyle{plain}
\bibliography{biblio}

\end{document}